\documentclass{amsart}

\usepackage{amsfonts,amssymb,amsmath,amscd,amsthm}
\usepackage[mathscr]{eucal}
\usepackage{mathrsfs}

\usepackage{mathbbol}

\usepackage{oldgerm,units}
\usepackage{wrapfig,epsfig}

\usepackage{ifthen}

\newtheorem{theorem}{Theorem}[section]
\newtheorem{proposition}[theorem]{Proposition}
\newtheorem{definition}[theorem]{Definition}
\newtheorem{algorithm}[theorem]{Algorithm}
\newtheorem{claim}[theorem]{Claim}
\newtheorem{lemma}[theorem]{Lemma}

\newtheorem{corollary}[theorem]{Corollary}

\newtheorem{example}[theorem]{Example}
\newtheorem{remark}[theorem]{Remark}

\newtheorem{observation}[theorem]{Observation}

\newtheorem{note}[theorem]{Note}

\newcommand{\rTheorem}[2]{\thmSpc \noindent \emph{\textbf{Theorem #1:} #2}\thmSpc}

\newcommand{\Ref}[1]{(\ref{#1})}

\newcommand{\Real}{\mathbb R}

\newcommand{\Net}{\mathbb N}

\newcommand{\Fld}{\mathbb K}

\newcommand{\graph}{\Gamma}

\newcommand{\Trop}{\mathbb T}
\newcommand{\RealInf}{{\Real \cup \{ \tUniS \}}}

\newcommand{\tUnit}{\mathbb U}
\newcommand{\etUnit}{\bar{\tUnit}}
\newcommand{\eReal}{\bar{\Real}}

\newcommand{\got}[1]{\frak{#1}}

\newcommand{\trop}[1]{\mathcal{#1}}

\newcommand{\nAmb}{\trop{N}}

\newcommand{\tC}{\trop{C}}
\newcommand{\tD}{\trop{D}}

\newcommand{\tF}{\trop{F}}

\newcommand{\tI}{\trop{I}}

\newcommand{\tP}{\trop{P}}

\newcommand{\tS}{\trop{S}}

\newcommand{\tZ}{\trop{Z}}

\newcommand{\tpF}{f}
\newcommand{\tpG}{g}
\newcommand{\tpH}{h}
\newcommand{\tpP}{p}
\newcommand{\tpQ}{q}
\newcommand{\tpFr}{\tpF_{\operatorname{r}}}
\newcommand{\tpFu}{\tpF_{\operatorname{u}}}
\newcommand{\tpGr}{\tpG_{\operatorname{r}}}
\newcommand{\tpGu}{\tpG_{\operatorname{u}}}

\newcommand{\tfF}{\tilde{f}}
\newcommand{\tfG}{\tilde{g}}

\newcommand{\tpFF}{F}
\newcommand{\tpGG}{G}

\newcommand{\To}{\longrightarrow }

\newcommand{\IFF}{\thSpc \Longleftrightarrow \thSpc}

\newcommand{\FigWidth}{2}

\newcommand{\oFace}{\tD}
\newcommand{\cFace}{\overline{\oFace}}

\newcommand{\subS}{\sqsubseteq}
\newcommand{\capS}{\sqcap}

\newcommand{\tUniS}{-\infty}

\newcommand{\uuu}[1]{#1^\nu}

\newcommand{\epiToMaxPlus}{\pi}

\newcommand{\polyEpiToMaxPlus}{\epiToMaxPlus_{*}}

\newcommand{\idl}[1]{\langle {#1} \rangle}

\newcommand{\al}{\alpha}
\newcommand{\bt}{\beta}
\newcommand{\gm}{\gamma}

\newcommand{\nXx}{{\bf x}}

\newcommand{\nTx}{{x_1,\dots,x_n}}

\newcommand{\nTz}{{z_1,\dots,z_n}}

\newcommand{\Val}{\operatorname{Val}}
\newcommand{\nVal}{Val}

\newcommand{\Deg}{Deg}

\def\indTrop{\operatorname{tr}}

\newcommand{\conC}{D}

\newcommand{\coVar}{\tC_{\indTrop}}
\newcommand{\uCoVar}{\widetilde{\coVar}}

\newcommand{\Id}{I}
\newcommand{\tId}{\Id_{\indTrop}}
\newcommand{\tIdA}{\got{a}}
\newcommand{\tIdB}{\got{b}}

\newcommand{\rtIdA}{\sqrt{\tIdA}}
\newcommand{\ptId}{\got{p}}
\newcommand{\mtId}{\got{m}}

\newcommand{\TrS}{\oplus}
\newcommand{\TrP}{\odot}

\newcommand{\TropSR}{(\Trop,\TrS,\TrP)}

\newcommand{\maxPlusAlg}{(\Real,max,+)}
\newcommand{\eMaxPlusAlg}{(\eReal,max,+)}

\newcommand{\CH}{\mathcal{CH}}
\newcommand{\uCH}{\overline{\CH}}

\newcommand{\OP}{\left(}
\newcommand{\CP}{\right)}

\def \aaa{a}
\def \aab{b}
\def \aac{c}

    \ifx\proof\undefined
    \newenvironment{proof}{
    \smallskip
    \noindent\emph{Proof.}}{\hfill\(\Box\)
    \bigskip
    } \fi

\newcommand{\bfem}[1]{\textbf{\emph{#1}}}

\newcommand{\ifdef}[3]{\ifthenelse{\equal{#1}{true}}{#2}{#3}}

\newcommand{\thSpc}{\; \; \;}

\newcommand {\secSpc} {\vskip 0cm}
\newcommand {\subSecSpc} {\vskip 0cm}

\newcommand {\parSpc} {\vskip 0.1cm}
\newcommand {\thmSpc} {\vskip 0.1cm}

\usepackage{epic}
\usepackage{eepic}

\def\pSkip{\vskip 1.5mm \noindent}
\def\piSkip{\vskip 1.5mm }

\def\bound{\partial}
\def\ru{(r,u)}

\def\hW{\widehat{W}}
\def\supp{\operatorname{Supp}}

\def\afSet{\tZ_{\indTrop}}
\def\cfSet{\tC_{\indTrop}}
\def\aSet{Z}
\def\cSet{C}

\def\Idl{\tI_{\indTrop}}

\def\bfa{{\bf a}}
\def\bfb{{\bf b}}
\def\bfc{{\bf c}}

\def\bfi{{\bf i}}
\def\bfj{{\bf j}}
\def\bfk{{\bf k}}
\def\bfx{{\bf x}}

\def\tnf{f^{\operatorname{t}}}
\def\gf{f^{\operatorname{g}}}

\def\ief{f^{\operatorname{i}}}
\def\ef{f^{\operatorname{e}}}
\def\eg{g^{\operatorname{e}}}
\def\ra{\text{a}}

\def\rad{rad}
\def\ldeg{\underline{\deg}}
\def\ttUnit{\tilde \tUnit}
\def\tTrop{\tilde \Trop}
\def\tReal{\tilde \Real}

\def\tf{\tilde f}
\def\tg{\tilde g}
\def\th{\tilde h}
\begin{document}


\title[Tropical Algebraic Sets, Ideals  and An Algebraic Nullstellensatz]
{Tropical Algebraic Sets, Ideals  and  An Algebraic
Nullstellensatz}


\author{Zur Izhakian}\thanks{The author has been supported by the grant of the High Council for
Scientific and Technological Cooperation between Israel and
France.}
\address{Department of Mathematics, Bar-Ilan University, Ramat-Gan 52900,
Israel} \email{zzur@post.tau.ac.il} \email{zzur@math.biu.ac.il}


\date{December 2007}


\keywords{ Idempotent semiring, max-plus algebra, tropical
algebraic geometry, algebraic sets and com-sets, polynomial
ideals, Nullstellensatz.}


\begin{abstract}
This paper introduces the foundations of the polynomial algebra
and basic structures for algebraic geometry over the extended
tropical semiring. Our development, which includes the tropical
version for the fundamental theorem of algebra, leads  to the
reduced polynomial semiring -- a structure that provides a basis
for developing a tropical analogue to the classical theory of
commutative algebra. The use of the new notion of tropical
algebraic com-sets, built upon the complements of tropical
algebraic sets, eventually yields the tropical algebraic
Nullstellensatz.
\end{abstract}

\maketitle


\tableofcontents

\section*{Introduction}\label{sec:Introduction}

The notion of tropical mathematics was introduced only in the past
decade \cite{Gelfand94,viro-2001}. Since then this theory has
developed rapidly and led to many applications
\cite{Einsiedler8311,Itenberg03,kontsevich-2000,MikhalkinEnumerative,Shustin2005,Speyer4218}.
A survey can be found in~\cite{Litvinov2005}. Tropical mathematics
is the mathematics over idempotent semirings,  the
\emph{\textbf{tropical semiring}} is usually taken to be $(\Real
\cup \{ \tUniS\}, \max, + \;)$; the real numbers, together with
the formal element $\tUniS$, equipped by the operations of maximum
and summation -- addition and multiplication respectively
\cite{kolokoltsov-2000}. The basic formalism of  tropical geometry
and been presented by Mikhalkin \cite{MikhalkinBook}.

\parSpc
The main goal of this paper is the development of another approach
to the basics of tropical polynomial algebra with a view to
tropical algebraic geometry, which is built on the
\emph{\textbf{extended tropical semiring}}, $\TropSR$, as has been
presented in \cite{zur05TropicalAlgebra}.
 This extension is obtained by taking two copies of the reals,
$$\eReal = \RealInf \quad \text{ and } \quad \etUnit = \uuu{\Real} \cup \{ \tUniS
\},$$ each is enlarged by $\{\tUniS\}$, and gluing them along
$\tUniS$ to define the set $\Trop = \eReal \cup \etUnit$. The
correspondence $\nu : \Real \to \tUnit$ is the identity map, so we
denote the image of $a \in \Real$ by $\uuu{a}$. Accordingly,
elements of $\tUnit$, which is  called the \bfem{ghost} part of
$\Trop$, are denoted as $\uuu{\aaa}$; $\Real$ is called the
\bfem{tangible} (or the \bfem{real}) part of $\Trop$. The map
$\nu$ is sometimes extended to whole $\Trop$,
\begin{equation}\label{eq:ghostMap}
    \nu: \Trop \To \etUnit,
\end{equation}
by declaring $\nu: \aaa^\nu \mapsto \aaa^\nu$ and $\nu : \tUniS
\mapsto \tUniS$; this map is called the \bfem{ghost map}.

The set $\Trop$ is then provided with the following total order
extending the usual order on $\Real$:
\piSkip
 $(i)$ $\tUniS \prec \al,$ $\forall \al \in \Trop$;
\piSkip
 $(ii)$ for any real numbers $\aaa < \aab$, we have $\aaa \prec
\aab,$ $\aaa \prec \uuu{\aab}$ and $\uuu{\aaa} \prec \aab$,
$\uuu{\aaa} \prec \uuu{\aab}$;
\piSkip
 $(iii)$ $\aaa \prec \uuu{a}$ for all $\aaa \in \Real$.
\pSkip
 (We use the generic notation
 $\aaa, \aab\in \Real$ and $\al,\bt \in \Trop$.)
  Then $\Trop$ is endowed with the two operations $\TrS$
and $\TrP$ , defined as follows: \pSkip
\qquad $
\begin{array}{cll}
  \al \TrS \bt & = &\left\{ \begin{array}{ll}
  max_{(\prec)}\{\al, \bt\}, & \al \neq \bt, \\[2mm]
  \uuu{\al}, & \al = \bt \neq \tUniS,
\end{array}\right. \\[2mm] \tUniS \TrS \tUniS & = & \tUniS,
\\[2mm] \aaa \TrP \aab & = & \aaa + \aab,\\[2mm]
a^\nu \TrP b & = & a \TrP b ^\nu = a^\nu \TrP b^\nu = (a+b)^\nu, \\[2mm]
 (\tUniS) \TrP \al  & = & \al \TrP (\tUniS) = \tUniS.
\end{array}
$ \pSkip

The semiring $(\Trop, \TrS, \TrP)$ modifies the classical
\bfem{max-plus algebra} and as has been proven, its arithmetic is
commutative, associative, and distributive. Note that while the
standard tropical semiring  $(\Real \cup \{ \tUniS\}, \max, + )$
is an idempotent semiring, since $\aaa \TrS \aaa = \uuu{\aaa}$,
the semiring $(\Trop, \TrS, \TrP)$ is not an idempotent semiring.
(The topology of $(\Trop, \TrS, \TrP)$ is more complected than the
Euclidean topology which is used on the standard tropical
semiring, the details are brought below in Section
\ref{sec:topolgy}.)

The connection with the standard tropical semiring is established
by the natural epimorphism of semirings,
\begin{equation}\label{eq:epiTopicalSemiRings}
    \epiToMaxPlus: \TropSR \; \To \;  (\Real \cup \{
\tUniS\}, \max, + \;), \\
 \end{equation}
given by  $\epiToMaxPlus:  \uuu{\aaa} \mapsto  \aaa$,
$\epiToMaxPlus: \aaa \mapsto \aaa$ for all $\aaa \in \Real$, and
$\epiToMaxPlus: \tUniS \mapsto \tUniS$. (We write
$\epiToMaxPlus(\al)$ for the image of $\al \in \Trop$ in
$\eReal$.) This epimorphism induces epimorphisms
$\polyEpiToMaxPlus$ of polynomial semirings, Laurent polynomial
semirings, and tropical matrices.

The fact that $(\Real, \TrP)$ is a group and $(\etUnit,
\TrS,\TrP)$ is an ideal provides $\Trop$ with a structure to which
much of the theory of commutative algebra (including polynomials
and determinants) can be transferred, leading to applications in
combinatorics, polynomials, Newton polytopes, algebraic geometry,
and convex geometry.

We start our discussion by observing the difference between
tropical polynomials and tropical polynomial functions, and study
the relation, which is not one-to-one correspondence, between
polynomials and functions. To overcome this miss-correspondence,
we determine the reduced polynomial semiring
$\tTrop[x_1,\dots,x_n]$ which is well behaved and allows an
analogous development of polynomial theory to that of the
classical case. This study includes polynomial factorizations and,
by introducing the tropical algebraic set
$$\afSet(\tpF)  = \{ \bfa \in \Trop^{(n)} \; | \; \tpF(\bfa) \in
\etUnit\}, \qquad \tpF \in \Trop[x_1,\dots,x_n],$$ one of our main
results is the \bfem{fundamental theorem of the
 tropical algebra} -- a tropical version that is similar
to the classical theorem.
\rTheorem{\ref{thm:FundamentalTheorem}}{The tropical semiring
$\Trop$ is algebraically closed (in tropical sense), that is,
$\afSet(\tpF) \neq \emptyset$ for any nonconstant $\tpF \in
\Trop[\nTx]$.}

The new notion of tropical com-set, defined as
$$ \cfSet(\tIdA) = \{ \conC_\tpF \; | \;
  \conC_\tpF \text{ is a  connected  component  of }
  \afSet(\tpF)^c \text{ of }  \tpF \in \tIdA \},$$
which are built upon the complements,  $\afSet(\tpF)^c$, of
tropical algebraic set $\afSet(\tpF)^c$, is central in our
development. The relation between com-sets and tropical ideals, is
the focal point for the
 tropical Nullstellensatz:
\rTheorem{\ref{thm:WeakTropicalNullstellensatz}}{(\textbf{Weak
Nullstellensatz}) Let $\tIdA \subset \tTrop[\nTx]$ be a finitely
generated proper ideal, then $\afSet(\tIdA) \neq \emptyset$.
Conversely, if $\afSet(\tIdA) = \emptyset$, then $\tIdA =
\tTrop[\nTx]$.}
\rTheorem{\ref{thm:AlgebraicTropicalNullstellensatz}}{(\textbf{Algebraic
Nullstellensatz})  Let $\tIdA \subset \tTrop[\nTx]$, where
$\ttUnit[x_1,\dots,x_n] \subseteq \tIdA$, be a finitely generated
tropical ideal, then
$\rtIdA = \tId(\cfSet(\tIdA)).$}
\noindent A similar context of the issues appear in this paper has
been raised in \cite[Qu. : A.16, C.2.A]{AIMworkshop}  and in
\cite[ Qu. 14]{AtlantaJointMeetings}.

\noindent \bfem{Notations}: In this paper we sometimes refer to
the standard arithmetic operations. To distinguish these
operations, the standard addition and the multiplication are
signed by $+$ and  $\cdot$ respectively. For short, we write $\aaa
\aab$ for  $\aaa \TrP \aab$.

\noindent  \textbf{\emph{Acknowledgement}}: The author would like
to thank \emph{Professor Eugenii Shustin} for his invaluable help.
I'm deeply grateful to him for his support and the fertile
discussions we had.



\secSpc
\section{The Topology of $\Trop^{(n)}$}\label{sec:topolgy}
Introducing a topology for $\Trop^{(n)}$, obtained as the product
topology on $\Trop$, in which the semiring's operations satisfy
continuity is essential for our future development.
Our topological setting is motivated by the following argument:
given a point $a \in \Real$ with a small neighborhood $\hW \subset
\Trop$, $a \in \hW$, pick $b \in \hW \cap \Real$, and consider the
sum $a \TrS b$ when $b \to a$. Then, in order to preserve the
continuity of $\TrS$ \,, $\hW$ must contain also the corresponding
ghost element $\uuu{a} \in \tUnit$. Later, we also want our
tropical sets to be closed sets.

Our auxiliary topology on the enlarged ghost part, $\etUnit =
\Real^\nu \cup \{ \tUniS \}$, is the Euclidean topology of the
half line $[0,\infty)$ in which $\TrS$ and $\TrP$ are continuous,
and closed sets are defined. The tangible part is concerned also
as having the Euclidean topology, but here, the topology is
partial, since $\TrS$ is continuous only for different elements.

Given a subset $U \subset \eReal$, we write $U^\nu $ for the the
corresponding  ghost subset $ \{ \uuu{u} \; | \; u \in U \}
\subset \etUnit$, recall that we identify $(\tUniS)^\nu$ with
$\tUniS$.
\begin{definition}\label{def:topology} A subset $\hW \subset \Trop$ is defined to be
\textbf{closed set} if \ $\hW = U \cup V^\nu$, where $U \subseteq
\eReal$ and $V^\nu \subseteq \etUnit$ satisfy:
\piSkip (i) $U^\nu ,V^\nu  \subseteq \etUnit$ are both closed sets
and,
\piSkip (ii) $U^\nu \subseteq V^\nu$.\pSkip
A set $W \subset \Trop$ is said to be \textbf{open} if its
complement is closed.
\end{definition}
\noindent (In particular, a closed  set may consist only of ghost
points, but when it includes a tangible point it must also contain
its ghost. Conversely, an open set can be pure tangible subset of
$\Real$.)

Using the decomposition $\hW = U \cup V^\nu$, it easy to verify
that finite unions and arbitrary intersections of the closed sets
are also of this form, accordingly, these sets form the closed
sets of our topology. Like in the standard case: the
\bfem{closure} of a set $W$ is the smallest closed set $\hW$
containing $W$, \bfem{connected set} $W$ is a set which cannot be
partitioned into two nonempty subsets such that each subset has no
points in common with the set closure of the other.

\begin{example} $ $ \pSkip
(i)  $\{ a, \uuu{a} \}$ and $\{ \tUniS\}$ are closed sets;  \pSkip
(ii) $\{ 1 \prec \al \prec 2 \; | \; \al \in \Real \}  $ is open
set; \pSkip
(iii)  $\{ \al \; | \; a \succ \al  \prec a^\nu \}$, for some $a
\in \Real$, is open set; \pSkip
(iv) $\{ 0 \preceq \al \preceq 1 \; | \; \al \in \Real \}  \cup \{
0^\nu \preceq \al \preceq 1^\nu \; | \; \al \in \tUnit \} $ is
closed set; \pSkip
\end{example}

 As mentioned earlier, having a topology on $\Trop$,  we
define the topology on $\Trop^{(n)}$ to be the product topology of
$\Trop$.

\secSpc
\section{Polynomials and Functions}

\subSecSpc
\subsection{The tropical polynomial semiring}\label{sec:TropVarieties}

The tropical semiring, $(\Trop[x], \TrS, \TrP)$, of polynomials in
one variable is defined to be all formal sums $f = \bigoplus _{i
\in \Net } \al _i x ^i$, with $\al_i \in \Trop$, for which almost
all $\al _i = -\infty,$ where we define polynomial addition and
multiplication in the usual way:
$$\left(\bigoplus _i \al _i x
^i\right)\left(\bigoplus_j \bt _j x ^j\right) = \bigoplus _k
\left(\bigoplus _{i+j=k} \al _i \bt _{k-j}\right) x ^k.$$
Accordingly, we write a polynomial $\bigoplus \al _i x ^i$ as
$\bigoplus_{i=0}^t \al _i x ^i$, when $\al _i = -\infty$ for all
$i>t$, and define its \bfem{degree} to be $t$. A term $\al_i x ^i$
is said to be a monomial of $f$ when  $\al_i \neq \tUniS$.

We sometimes write $x ^\nu$ for $0^\nu x.$ Note that, since $0$ is
the multiplicative unit of $\Trop$, we write $x^i$ for $0 x ^i$
and say a polynomial is \bfem{monic} if its coefficient of highest
degree, which we call the leading coefficient, is $0$. We identify
$\al x^0$ with $\al$, for each $\al \in \Trop$; thus we may view
$\Trop \subset \Trop[x]$. The elements of $\tUnit[x]$ form the
ghost part of $\Trop[x]$ where $\Real[x]$ is its tangible part;
that is, polynomials of each part have, respectively, only ghost
coefficients or only tangible coefficients.

The polynomial semiring $\Trop[x_1,\dots,x_n]$ is defined
inductively, as $\Trop[x _1, \dots, x _{n-1}][x _n]$; a typical
polynomial, as usual, is
$$f = \bigoplus \al_{i_1, \dots, i_n} x _1^{i_1}\cdots  x _n^{i_n}.$$
We write $\bfi$ for multi index $(i_1, \dots, i_n)$ and let $\bfx
= (x_1,\dots,x_n)$, and thus write $f = \bigoplus \al_\bfi
\bfx^\bfi$,  for short. The \bfem{support} of a polynomial $f$ is
define to be those $\bfi$ for which $\al_\bfi \ne \tUniS$, that is
$$\supp (f) = \{ \; \bfi \; | \; \al_{\bfi}  \neq \tUniS \}.  $$

Given $f \in \Real[x_1,\dots,x_n]$, i.e. $\al_\bfi = a_\bfi \in
\Real$ for all $\bfi$, the corresponding polynomial $f = \bigoplus
a^\nu_\bfi \bfx^\bfi$ is denoted as $f^\nu$. Moreover, $f$ can be
decomposed uniquely according to its tangible part $\tnf$ and its
ghost part $\gf$, and written uniquely as $f = \tnf \TrS \gf$. We
call this \bfem{(t,g)-decomposition} of $f$, clearly this
decomposition is unique. If $f = \tnf$ then $f$ is said to be
\bfem{tangible polynomial}, and is said to be \bfem{ghost
polynomial} when $f = \gf$.

\begin{remark}
The tangible part $\Real[x_1,\dots, x_n]$  is not closed under the
semiring  operations. Moreover, there are $f \in \Real[x_1,\dots,
x_n]$ for which $f^k \notin \Real[x_1,\dots, x_n]$ for some
positive $k \in \Net$; for example take $f = x \TrS 1$ then $f^2 =
x^2 \TrS 1^\nu x \TrS 2$ which is a non-tangible polynomial (but
is not ghost polynomial).

A power of a non-ghost polynomial, i.e. it has a tangible
monomial, can be a ghost polynomial; for  example take $ \tpF =
\uuu{0} x^2 \TrS 1x \TrS \uuu{2},$ then
$$
  \tpF^2 =  (\uuu{0} x^2  \TrS 1x \TrS \uuu{2})^2 =  \uuu{0}( x^4   \TrS 1 x^3 \TrS
           2x^2  \TrS 3x \TrS 4),
 $$
 which is a ghost polynomial. On the other hand, $\tUnit[x_1,\dots, x_n]$ is
  closed under addition
 and under the
multiplication with any element of $\Trop[\nTx]$; therefore, as
will be seen later, is a semiring  ideal.
\end{remark}
\noindent We note that whenever $fg = \tUniS$ then $f= \tUniS$ or
$g = \tUniS$, and thus  the only element $\tpF \in \Trop[\nTx]$
for which $\tpF^k = \tUniS$ is $\tUniS$ itself.

Recall that $\Trop$ lacks subtraction, and therefore we don't have
cancelation of monomials;  this property is expressed in degree
computations that always satisfy the  rules:
$$
  \deg(\tpF \tpG)  =  \deg(\tpF) + \deg(\tpG) \quad \text{ and }
  \quad
  \deg(\tpF \TrS \tpG)  =  \max\{ \deg(\tpF), \deg(\tpG)\}.
$$
(This is different from the classical theory in which $\Deg(\tpF +
\tpG) \leq max\{ \deg(\tpF), \deg(\tpG)\}$, in the tropical case
we always have equality.) For the same reason, for a given $\tpF
\in \Trop[\nTx]$ one can define the  ``lower degree'' to be
$$\ldeg(\tpF) = min\{ \deg(h) \; | \; h \; is \; a \; monomial \; of \; \tpF \}.$$
which then satisfies
$$
  \ldeg(\tpF \tpG) =  \ldeg(\tpF) + \ldeg(\tpG) \quad \text{ and }
  \quad
  \ldeg(\tpF \TrS \tpG)  =  \min\{ \ldeg(\tpF), \ldeg(\tpG)\}.
$$
Clearly, we always have $\ldeg(\tpF) \leq \deg(\tpF)$ and both
can only increase by preforming operations over polynomials.

A \bfem{tropical homomorphism} of tropical polynomial semirings
$$\varphi: (\Trop[x_1,\dots,x_n],\TrS,\TrP) \ \To \  (\Trop[x_1,\dots,x_m],\TrS,\TrP)$$ is a semiring  homomorphism
 $\varphi : \Trop[x_1,\dots,x_n] \setminus \{ \tUniS \} \to\Trop[x_1,\dots,x_m] \setminus \{ \tUniS \}$
such that $\varphi(f^\nu) = (\varphi(f))^\nu$ for any $f \in
\Trop[x_1,\dots,x_n]$; accordingly, we have $\varphi
(\tUnit[x_1,\dots,x_n]) \subset \tUnit[x_1,\dots,x_m] $. (Within
this definition we include the case of $m =0$, that is $\varphi:
(\Trop[x_1,\dots,x_n],\TrS,\TrP)  \to \Trop$.)  The \bfem{tropical
kernel}, $\ker \varphi$, is the preimage
 of $\etUnit[x_1,\dots,x_m]$.
We call $\varphi$ a \bfem{ghost injection} if $\ker \varphi =
\etUnit[x_1,\dots,x_n]$ and say that $\varphi$ is a \bfem{tropical
injection} if $\varphi$ is 1:1 and
 is a ghost injection.

Given a point  $\bfc = (c_1, \dots, c_n) \in \Trop^{(n)}$, there
is a tropical homomorphism $\varphi_\bfc:
 \Trop[x_1,\dots, x_n]\to \Trop,$ given by sending
 $$\varphi_\bfc:\bigoplus _i \al_{i_1, \dots, i_n} {x  _1}^{i_1}\cdots  {x _n}^{i_n}
 \  \longmapsto \
 \bigoplus _i \al_{i_1, \dots, i_n} {c _1}^{i_1}\cdots  {c _n}^{i_n},$$
which we call the \bfem{substitution homomorphism} (with respect
to $\bfc$). (Note that in tropical algebra, as usual,
${c_j}^{i_j}$ means the tropical product of $c_j$ taken $i_j$
times, which is just $i_j \cdot c_j$ in the classical notation.)
We write $f(\bfc)$ for the image of the polynomial $f$ under
substitution to $\bfc$.

In our philosophy, elements of $\etUnit$ are those to be ignored,
accordingly we define a root of polynomial, in the tropical sense:
\begin{definition}
An element $\bfa \in \Trop^{(n)}$  is a \textbf{root} of $f$ if
$f(\bfa) \in \etUnit$, i.e if $f \in \ker \varphi_{\bfa }$, where
$\varphi_{\bfa }$ is the tropical substitution homomorphism
$\varphi_{\bfa }: (x_1, \dots, x _n) \mapsto \bfa$.
\end{definition} \noindent (Note that we include $\tUniS$ as a proper root
since later we want to study connectedness of sets of roots and
their complements which coincides with our topological setting.)
By this definition, given a ghost polynomial $f = \gf$, then any
$\bfa \in \Trop^{(n)}$ is a root of $f$. Note the we can also have
non-ghost polynomials for which any $\bfa \in \Trop^{(n)}$ is a
root (take for example $\tpF = \uuu{0}x^2 \TrS-1x \TrS \uuu{0}$).
However, our main interest is in non-ghost polynomials, mainly in
tangible ones.

\begin{remark}
Suppose $f = \tnf \TrS \gf$ is the (t,g)-decomposition of a
polynomial $f$ into tangible and ghost parts. Then any root $\bfa$
of $\tnf$ is a root of $f$. Indeed, $f(\bfa) = \tnf (\bfa) \TrS
\gf(\bfa)$, and each part is in $\etUnit$.
\end{remark}

\begin{lemma}\label{lem:cl1} For
any nonconstant polynomial $f \in \Trop[x]$ without a constant
monomial and for any $\ra^\nu \ne \tUniS$ in $\tUnit$, there
exists $ r \in \Trop$ with $f(r) =\ra^\nu \in \tUnit$,
\end{lemma}
\begin{proof}  Write $f = \bigoplus \al_i x^i.$ For each $i>0$,
there is some $r_i \in \Real$ such that $ \al_i ({r_i}^i)^\nu =
\ra^\nu$. Indeed, assume $\al_i$ is tangible then using the ghost
map \Ref{eq:ghostMap} (written in the standard arithmetic) we have
$$\nu (r_i) = \frac 1i(\ra - \al_i)^\nu.$$
 Now, take $r$ among these $r_i$ such that $\nu (r)$ is minimal
among $r_i^\nu, \ 1 \le i \le t.$ Then $ (f(r))^\nu = \ra ^\nu .$
\end{proof}

 Using Lemma \ref{lem:cl1}, we can state the
\bfem{Fundamental theorem of the  tropical algebra} -- a tropical
version that is similar to the classical theorem.
\begin{theorem}\label{thm:FundamentalTheorem}
The tropical semiring $\Trop$ is algebraically closed in tropical
sense, that is any nonconstant tropical polynomial $\tpF \in
\Trop[\nTx]$ has a root.
\end{theorem}

\begin{proof} Assume $f \in \Trop[x]$, if $f$ has a single nonconstant monomial $h_i = \al_ix
^i$, or $f = x^i g$, for some $i > 0$, then $h(\ra^\nu) \in
\tUnit$ and  respectively $(\ra^\nu)^i g(\ra^\nu) \in \tUnit$, for
any $\ra^\nu$, and we are done. Suppose $f = \bigoplus _{i=0}^m
\al _{i} x ^i,$ we may assume that $\al_0 \ne -\infty.$
 Write $f = g
\TrS \al _0.$ If $\al_0 \in \tUnit $, then for any ghost
$\ra^\nu,$ $f(\ra^\nu) \in \tUnit,$ so $\ra^\nu$ is a root. Thus,
we may assume that $\al_0 \notin \etUnit$. By lemma \ref{lem:cl1},
there is some $r$ such that $\nu (g (r)) = \al_0^\nu,$ which
implies $ f(r)= \al_0^\nu \TrS \al _0 \in  \tUnit.$ The
generalization to $f \in \Trop[\nTx]$ is clear.
\end{proof}

\begin{remark} In the familiar tropical semiring $(\Real,\max,+)$ roots are
not defined directly and are realized  as the points on which the
evaluation of a polynomial is attained by at least two of its
monomials. In other word, the roots are simply the domain of
non-differentiability of the corresponding function.
Unfortunately, using this notion, $(\Real,\max,+)$ is not
algebraically closed in the tropical sense; take for example a
polynomial having a single monomial.
\end{remark}

\subSecSpc
\subsection{Tropical polynomial functions}
As mentioned earlier, in the tropical world the correspondence
between polynomials and polynomial functions is not one-to-one,
mainly due to convexity matters, and a function can have many
polynomial descriptions; for example consider the family of the
polynomials
$$ f_t = x^2 \TrS t x \TrS 0,$$
where $t \leq 0 $ serves as a parameter, all the members of this
family describe the same function. We denote by $\psi_f$ the
tropical function corresponding to a polynomial $f \in
\Trop[x_1,\dots,x_n]$,  that is $ \psi_f : \bfa \mapsto f(\bfa)$
written as $\psi_f(\bfa)$. We denote by $\tF(\Trop^{(n)})$ the
semiring of polynomial functions
$$ \tF(\Trop^{(n)}) = \{ \psi_f : \Trop^{(n)} \To \Trop \ | \ f \in
\Trop[x_1,\dots,x_n] \}.$$ (The operations
$$(\psi_f \TrP \psi_g)(\bfa)  = \psi_f(\bfa) \TrP \psi_g(\bfa),
  \qquad and \qquad
  (\psi_f \TrS \psi_g)(\bfa)  = \psi_f(\bfa) \TrS \psi_g(\bfa),$$
of $\tF(\Trop^{(n)})$  are defined point-wise.)
 Given a tropical function, the central
idea for the further development is finding the best
representative among all of its polynomials descriptions.

\begin{definition}\label{def:R-equivalent}
 Two polynomials $f,g\in \Trop[x_1,\dots,x_n]$ are said to be \textbf{equivalent},
 denoted as $f \sim g$,  if they take on the same values,
  that is $f(\bfa) = g(\bfa)$, for any $\bfa = (\ra_1,\dots,\ra_n) \in \Trop^{(n)}$
  (i.e. in function view $\psi_f = \psi_g$). \end{definition}

\begin{example}\label{exmp:sim} For all $a,b \in
\Real$, $a \neq b$, the following relations hold true: \pSkip (i)
$x \TrS a \nsim x \TrS \uuu{a}$,  \pSkip (ii) $x \TrS a \nsim x
\TrS b, $ \pSkip (iii) $(x \TrS \al)^2 \sim x^2 \TrS \al^2$.
\pSkip

To prove (iii), write  $f = (x \TrS \al)^2 = x^2 \TrS \al^\nu x
\TrS a^2$ and $g = x^2 \TrS \al^2$. Suppose $f \nsim g$, this
means that there is some $\ra \in \Trop$ for which $f(\ra) \neq
g(\ra)$ and thus $f(\ra) = \al^\nu \ra \succ  \al^2 \TrS \ra^2$.
But, if $\ra \succ \al$ we have $\ra^2 \succ \al^\nu \ra$, and
when $\ra \prec \al$ we get $\al^2 \succ \al^\nu \ra$. which is a
contradiction. (In the case of $ \ra = \al$ we have $f(\al) =
g(\ra) = \al^\nu$.) Namely, $f(\ra) = \al^2 \TrS \ra^2$ for any
$\ra \in \Trop$.

\end{example}

Given $f \in \Trop[x_1,\dots.x_n]$, the graph $\graph(\tpF)$ of
$f$ is defined to be
$$\graph(\tpF) = \{ (\ra_1, \dots, \ra_n, \tpF(\bfa)) \ : \ \bfa = (\ra_1, \dots, \ra_n)
 \in \Trop^{(n)}  \}  \subset \Trop^{(n+1)}.$$
We write $\graph(\tpF|_\Pi  )$  for the restriction of to the
subdomain  $\Pi \subset \Trop^{(n)} $, that is $$\graph(\tpF|_\Pi
) = \{ (\ra_1, \dots, \ra_n, \tpF(\bfa)) \ : \ \bfa  \in \Pi \}.$$
Accordingly, we have  the following relations
\begin{lemma}\label{lem:eq} Let $f \in \Trop[x_1,\dots.x_n]$ then
\begin{equation*}\label{eq:tropPolynomialFun}
\tpF \sim  \tpF' \IFF  \graph(\tpF) = \graph(\tpF')
\end{equation*}
and
$$  \deg(\tpF) = \deg(\tpF') \quad \text{and} \quad
  \ldeg(\tpF) =   \ldeg(\tpF').
$$
\end{lemma}

\begin{proof} The first relation is by definitions,
$\tpF \sim  \tpF'$ if and only if  $f(\bfa) = f'(\bfa)$ for all
$\bfa = (\ra_1,\dots,\ra_n) \in \Trop^{(n)}$ if and only if
$(\ra_1,\dots,\ra_n, f(\bfa)) = (\ra_1,\dots,\ra_n, f'(\bfa)) $
for all $\bfa = (\ra_1,\dots,\ra_n) \in \Trop^{(n)}$ if and only
if $\graph(\tpF) = \graph(\tpF')$.

Suppose $f \sim f'$ and $\deg(\tpF)
> \deg(\tpF')$.
Write $\bfi = (i_1,\dots,i_n)$, $\bfj = (j_1,\dots,j_n)$ for the
multi-indices  and let $\al_\bfi \bfx^\bfi $ and  $\al_\bfj
\bfx^\bfj $ be respectively the monomials of highest degree of $f$
and $f'$. Then $i_s > j_s$ for some $s = 1,\dots,n$, which implies
that for a point $\bfa \in \Trop^{(n)}$ whose $s$'th coordinate is
sufficiently large we have $f(\bfa) \succ f'(\bfa) $ -- a
contradiction. To prove that $\ldeg(\tpF) = \ldeg(\tpF')$, we use
the same argument only by considering the monomials of lowest
degree with respect to a point having a sufficiently small
coordinate.
\end{proof}

\begin{remark}
Assume $\tpF \in \Trop[x_1,\dots,x_n]$ is tangible (resp. ghost)
and $f \sim f'$, then $\tpF'$ needs not  be also tangible (resp.
ghost); for example  $x^2 \TrS \al^\nu x \TrS \al^2 \sim x^2 \TrS
\al^2$ (cf. Example \ref{exmp:sim}).
\end{remark}

Instead of polynomials, we are interested in their equivalence
classes. There is a natural representative for each equivalence
class. Given a polynomial $f = \bigoplus_\bfi \al_\bfi \bfx^\bfi$
having a monomial $h =  \al_\bfj \bfx^\bfj$, we denote by $f
\setminus h$ the polynomial   $\bigoplus_{\bfi \neq \bfj} \al_\bfi
\bfx^\bfi $.

\begin{definition}\label{def:essentialPart} A polynomial $f \in \Trop[x_1,\dots, x_n]$
\textbf{dominates} $g$ if $f(\bfa) \TrS g(\bfa) = f(\bfa)$ for all
$\bfa = (\ra_1,\dots,\ra_n) \in \Trop^{(n)}$, i.e. $\psi_f \TrS
\psi_g = \psi_f$.
 A monomial $h$ of $f$ dominated by $f \setminus h$ (which not empty)
is called \textbf{inessential}; otherwise $h$ is said to be
\textbf{essential}. The \textbf{essential part} $\ef$ of a
polynomial
 $f = \bigoplus \al _\bfi \bfx ^\bfi$ is the sum of those monomials
$\al_\bfi \bfx ^\bfi$ that are essential, while its inessential
part $\ief$ consists of the sum of all inessential $\al_\bfi \bfx
^\bfi$. When $f = \ef$, $f$ is said to be an \textbf{essential
polynomial}.
\end{definition}
For example, $x^2 \TrS 2$ is the essential part of  $x^2 \TrS 0x
\TrS 2$, where $0x$ is inessential monomial. In other words, the
essential part  consisting  of all the monomials which are need to
obtain a same polynomial function. Namely, from the function point
of view, to obtain $\ef$ we cancel out all the unnecessary
monomials of $f$.

\begin{lemma}\label{lem:fef}
For any  $f \in \Trop[x_1,\dots,x_n]$,  $f \sim \ef$.
\end{lemma}
\begin{proof} Let $f=\bigoplus_i h_i$ and assume that $f \nsim \ef$.
Then, there is some $\bfa \in \Trop^{(n)}$ for which $f(\bfa) \neq
\ef(\bfa)$. This means that $f \setminus h_i$ does not dominate
some monomial $h_i$ and this monomial is not part of $\ef$.
Namely, $f$ contains an essential monomial $h_i$ which is not in
$\ef$. This contradicts the construction of $\ef$.
\end{proof}

\begin{proposition}\label{pop:uniqe}
The essential part, $\ef$, of a polynomial $f$ is unique.
\end{proposition}
\begin{proof}   Assume that $f$ have two different essential parts,
say $\ef$ and ${\ef}'$, then $\graph(\ef) \neq \graph ({\ef}')$.
But then, by Lemma \ref{lem:fef}, $f \sim \ef$ and $f \sim
{\ef}'$, and by Lemma \ref{lem:eq},  $\graph(\ef) = \graph
({\ef}')$ -- a contradiction.
\end{proof}

Integrating Lemma \ref{lem:eq}, Lemma \ref{lem:fef}, and
Proposition \ref{pop:uniqe} we conclude:
\begin{corollary}\label{cor:essential}
 $f \sim g$  if and only if $\ef = \eg$.
\end{corollary}

Clearly, $\sim$ is an equivalence relation, so $\ef$  serves as a
canonical representative for the equivalence class $$\tC_f = \{f'
\in \Trop[x_1, \dots, x_n] : \ f' \sim f\}.$$  Thus, each
equivalence class under $\sim$ has a canonical (essential)
representative. One can use these representatives to establish the
one-to-one correspondence:
$$\Trop[x_1, \dots, x_n]/_\sim  \To \tF(\Trop^{(n)}).$$
Yet, we are looking for a better representative since these
representatives are not suitable for the purpose of factorization.

\begin{note}\label{obs:geometryOfGraphs}
Assume the essential part of $\tpF$ is tangible and is comprised
of $m$ tangible monomials (i.e. $\ef \in \Real[\nTx]$),
considering $\ef$ over $\eMaxPlusAlg$, then
$\graph_{\maxPlusAlg}(\ef) \subset \Real^{(n+1)}$ is a convex
polyhedron having $m$ faces $\cFace_i$ of codimension 1. On the
other hand, the tangible part of $\graph(\ef)$ over $\Real^{(n)}$
(i.e. $\graph(\ef) \cap \Real^{(n+1)}$) consists of the same faces
$\oFace_i$ as those of $\graph_{\maxPlusAlg}(\ef)$ but without
their boundaries (in the view of the Euclidean topology). These
boundaries ``pass'' to $\Real^{(n)} \times \tUnit$, so in
$\Real^{(n+1)}$ the faces $\oFace_i$ are open sets. In other
words, using the Euclidean topology for $\Real^{(n+1)}$,
$\graph_{\maxPlusAlg}(\ef)$ is the closure of $ \graph(\ef|_
{\Real^{(n)}})$. Note that this is true only for tangible $
 \ef$, yet we always have the onto projection
$
  \graph(\ef) \; {\To} \; \graph_{\maxPlusAlg}(\ef)
$ and for any $\Pi_i  \in \{\Real,\tUnit\} \times \dots \times
\{\Real,\tUnit\}$, we have the isomorphism
\begin{equation}\label{eq:homomorphismOfCopy}
  \graph(\ef|_{\Pi_i}) \; \overset{\sim}{\To} \; \graph_{\maxPlusAlg}(\ef).
\end{equation}
In general, over $\Trop^{(n)}$, we have $2^n$ subgraphs, $\graph
(\ef|_{\Pi_i})$, each is isomorphic to
$\graph_{\eMaxPlusAlg}(\ef)$ in $\eMaxPlusAlg$.
\end{note}

Suppose $f = \bigoplus_\bfi \al _\bfi \bfx^\bfi \in
\Trop[x_1,\dots,x_n]$, we identify each monomial $\al_\bfi
\bfx^{\bfi}$ (for $\bfi = (i_1, \dots, i_n)$) with the point
$$ (i_1, \dots, i_n, \epiToMaxPlus( \al _\bfi)) \in \Net
^{(n)} \times \Real \subset \Real^{(n+1)}.$$
Let $C_f$ be the polyhedron determined by the points
$$\{ (i_1, \dots, i_n, \epiToMaxPlus(\al _\bfi)): \bfi \in \supp (f)\},$$
which we call the \bfem{vertices} of $C_f$, and take the convex
hull $\CH_f$ of these vertices. We say that a vertex  is tangible
(resp. ghost) vertex if it corresponds to a tangible monomial
(resp. ghost monomial).

 Let $A_\bfj \subset \CH_f$ be the set of
points whose first $n$ coordinates are equal. The point $\bfa_\bfj
= (j_1,\dots,j_n, a) \in A_\bfi$ whose $(n+1)$'th coordinate is
maximal among all the points of $A_\bfj$ is said to be an
\bfem{upper point} of $\CH_f$. The upper part of $\CH_f$,
consisting  of all the upper points in $\CH_f$, is called the
\bfem{essential complex} of $f$ and is denoted
 $\uCH_f$. The points of  $\uCH_f$ of the form
$\{(i_1, \dots, i_n, \epiToMaxPlus(\al _\bfi)): \bfi \in \Net
^{n}\}$ are called \bfem{lattice points}. For example, when $f =
x^2 + 2,$ the lattice points are (2,0), (1,1), and (0,2).

Note that the essential complex can be consisted of both tangible
and ghost vertices, in particular the essentiality of a vertex
(and of monomial, as will be seen later) is independent on being
tangible or ghost.

In fact the structure described above can be understood in the
more winder context of the Newton polytope \cite{IMS}. Recall that
the Newton polytope $\Delta_f$, of $f = \bigoplus_{\bfi} \al_\bfi
\bfx^\bfi$ is the convex hull  of the $\bfi$'s in $\supp(f)$. By
taking the onto projection, which is obtained by deleting the last
coordinate, of the non-smooth part of $\uCH_f$ (that is a
polyhedral complex) on $\Delta_f$ the induced polyhedral
subdivision $S_f$ of $\Delta_f$ is obtained. Thereby, a dual
geometric object having combinatorial properties is produced. This
object plays a major  role in classical tropical theory and it
being used in many applications
\cite{Itenberg03,MikhalkinEnumerative,Shustin1278}.

\begin{lemma}\label{lem:essential} There is a one-to-one correspondence between the
vertices of the essential complex $\uCH_f$ of $f$ and the
essential monomials of $f$.
\end{lemma}
\noindent Vertices of the essential complex  $\uCH_f$ are in
one-to-one correspondence with the vertices of the induced
subdivision $S_f$ of Newton polytope $\Delta_f$. (The latter are
precisely the projections of the vertices of $\uCH_f$ on
$\Delta_f$.) The proof is then obtained by the one-to-one
correspondence between vertices of $S_f$ and essential monomials
of $f$ \cite{Shustin1278}.

Note that $\uCH_f$ may contain lattice points not corresponding to
monomials of the original polynomial $f$. For instance, take $f =
x^2 + 2,$ then the lattice point (1,1) does not correspond to a
monomial of $f$. In general, the inessential part of $f$ does not
appear in $\uCH_f$ as vertices but it may appear as points that
lie on its faces. A vertex of $\uCH$ is called \bfem{interior} if
its projection to $\Delta_f$ is not a vertex (but is still a
vertex of $S_f$).
 We say the monomial $h_\bfi =\al _\bfi \bfx
^\bfi$ is \bfem{quasi-essential} for $f$ if $(i_1, \dots, i_n,
\epiToMaxPlus(\al_\bfi))$ lies on $\uCH_f$ and is not a vertex.
This has the following interpretation:

\begin{lemma}
An inessential monomial is  quasi-essential if any (arbitrarily
small) increase of its coefficient makes it essential.
\end{lemma}

\begin{proof} Let $\al_i \bfx^\bfi$ be a quasi-essential monomial.
Any arbitrarily small increasing of its coefficient $\al_\bfi$
makes the corresponding lattice point $(i_1,\dots,i_n,
\epiToMaxPlus(\al _\bfi))$ of $\uCH$ a vertex. Then, by Lemma
\ref{lem:essential}, $\al_i \bfx^\bfi$ becomes essential.
\end{proof}

\begin{remark}\label{rem:convexity} Summarizing the above discussion, we see that  the polynomial corresponding to the upper
part of $\CH_f$ is precisely the essential part of $f$, and in
particular $\uCH_{\ef} = \uCH_f$. Thus, two polynomials are
equivalent iff they have the same essential part iff their
essential complexes, including their indicated tangible/ghost
vertices, are identical.
\end{remark}

Any  $\tpF \in \Trop[\nTx]$ can be written uniquely as $$\tpF =
\tpFr \TrS \tpFu $$ with $\tpFr,\tpFu \in \Real[\nTx]$, we call
this form the \bfem{$\textbf{\ru}$--decomposition} of $f$. To
obtain this decomposition, just take each ghost monomial $\al_\bfi
\bfx^\bfi$ (i.e. $\al_\bfi \in \tUnit$) and replace it by the two
tangible copies $\epiToMaxPlus(\al_\bfi) \bfx^\bfi$, i.e.
\begin{equation}\label{eq:ruDecom}
\al_\bfi \bfx^\bfi \; \rightsquigarrow \; \epiToMaxPlus(\al_\bfi)
\bfx^\bfi \TrS \epiToMaxPlus(\al_\bfi) \bfx^\bfi.
\end{equation}
Then, take one copy from each pair of these monomials to create
$\tpFu$, the remaining monomials are ascribed to $\tpFr$, in
particular $\tpFu = \polyEpiToMaxPlus(\gf)$ and $f = \tpFr$ if $f$
is tangible.  In this view we have the following:
\begin{proposition}
$f \sim g$ if and only if $\uCH_{\tpFr} = \uCH_{\tpGr}$ and
$\uCH_{\tpFu} = \uCH_{\tpGu}$.
\end{proposition}
\begin{proof}  By Corollary
\ref{cor:essential} $f \sim g$ iff $\ef = \eg$, so we may assume
$f$ and $g$ are essential.  Since $\ru$-decomposition is unique,
we get $\tpFr = \tpGr$ and $\tpFu = \tpGu$, where all $\tpFr$,
$\tpGr$, $\tpFu$, and $\tpGu$ are essentials. Thus, $\uCH_{\tpFr}
= \uCH_{\tpGr}$ and $\uCH_{\tpFu} = \uCH_{\tpGu}$.
\end{proof}

\subsection{The representatives of polynomial classes}

Next we want to identify the best canonical representative of a
class of equivalent polynomials. Note that we already have a
canonical representative, which is the common essential part of
all the class members. Yet, we are looking for a better
representative which, as will be seen later, is useful for easy
factorization; for this purpose we need the following:

\begin{definition}\label{def:fullPoly} A polynomial
 $f \in \Trop[x_1,\dots,x_n]$ is called \textbf{full} if every lattice
point lying on $\uCH_f$ corresponds to a monomial which is either
essential or quasi-essential, and furthermore, every nonessential
monomial is ghost; a full polynomial $f$ is \textbf{tangible-full}
if $\ef$ is tangible. The \textbf{full closure} $\tilde f$ of $f$
is the sum of $\ef$ with all the quasi-essential  monomials of $f$
taken ghost.
\end{definition}
\noindent By this definition, the full closure is unique, and
therefore $\tilde f$ is also canonical representative of $\tC_f$.
We call $\tilde f$ the \bfem{full representative} of $\tC_f$, this
representative plays a major role  in our future development.

\begin{remark} When $f$ is a polynomial consisting of  a single
monomial, then $f$ is (full) essential and we always have $f =
\tilde f$.

\end{remark}
 \begin{example}  $ $ \pSkip
(i) $x^2 \TrS 1^\nu x \TrS 0$ is tangible-full essential; \pSkip
(ii) $x^2 \TrS 1^\nu x \TrS 0^\nu$ is full essential; \pSkip (iii)
$x^2 \TrS 0^\nu x \TrS 0$ is full but not essential; \pSkip (iv)
$x^2 \TrS 0 x \TrS 0$ is not full since $0x$ is tangible; \pSkip
(v) $x^2 \TrS 0^\nu x \TrS 0$ is the full closure of $x^2 \TrS 0 x
\TrS 0$ and $x^2 \TrS 0$.
 \end{example}

By the construction of $\uCH_f$ and the fact that the full
polynomials contain all the monomials corresponding to lattice
points of their essential complexes we have the following:

\begin{lemma}\label{rem:desending} Any full polynomial $f \in
\Trop[x]$ (which is not a monomial) corresponds to a descending
sequence of tangible elements $m_1,\dots, m_t$, where $t = \deg f
- \ldeg f $, which is defined uniquely by the slopes of the series
of edges $e_1,\dots, e_t$ of $\uCH_f \subset \Real^{(2)}$, each
$e_i$ is determined by the pair $(i-1,\epiToMaxPlus(\al_{i-1}))$
and $(i,\epiToMaxPlus(\al_{i}))$.

\end{lemma} The descending sequence of tangible elements $m_1,\dots,
m_t$ is denoted by $M_f$. Note that $M_f$ is not necessarily
strictly descending and it might have identical adjacent elements.
The sequence of edges is denoted by $E_f$.
\begin{proof}
 Recall that since $f$
is full, it has exactly $t+1$ monomials, and by the construction
of $\uCH_f$ it also contains $t+1$ lattice points (not all of them
need to be vertices). The sequence $M_f$ is descending due to the
convexity of $\uCH_f$. Since otherwise, assume $m_{i+1} > m_{i}$,
for some $i = 1,\dots,t-1$, and observe the corresponding lattice
points
$$ (i-1,\epiToMaxPlus(\al_{i-1})  ),  \ (i, \epiToMaxPlus(\al_{i})), \ (i+1, \epiToMaxPlus(\al_{i+1})),$$
which by assumption should satisfy
$$ \epiToMaxPlus(\al_{i+1})- \epiToMaxPlus(\al_{i}) >   \epiToMaxPlus(\al_{i}) - \epiToMaxPlus(\al_{i-1}).$$
(Here use the standard notation to describe the slopes of the
edges since we work only on $\Real^{(2)}$.) But this means, due to
the convexity of $\uCH_f$, that $(i, \epiToMaxPlus(\al_{i}))
\notin \uCH_f$ and thus, is not a lattice point.
\end{proof}

\begin{remark} Clearly, the lemma holds true for the
full closure of any $f \in \Trop[x]$. Moreover,  one can state
Lemma \ref{rem:desending} for any essential polynomial $f \in
\Trop[x]$, but in this case the number of monomials of $f$ will be
less or equal to $t$.
\end{remark}

\begin{definition}\label{def:reduced} The \textbf{reduced domain} $\tTrop[x_1,\dots,x_n]$ of
$\Trop[x_1,\dots,x_n]$ is the set of  full elements, where
addition and multiplication are defined by taking the full
representative of the respective sum or product in
$\Trop[x_1,\dots,x_n].$ In other words,  we define
$$f \TrS g = \widetilde {f \TrS g}, \qquad f \TrP g = \widetilde {f \TrP
g},$$ for $f,g \in \tTrop[x_1,\dots,x_n] $, and sometimes call
$\tilde\Trop[x_1,\dots,x_n]$ the \textbf{reduced polynomial
semiring}. Accordingly,
 $\tilde\Real[x_1,\dots,x_n]$ is the set of tangible-full elements, and $\tilde \tUnit[x_1,\dots,x_n]$ is the set of full
elements, all of whose coefficients are ghosts.
\end{definition}

Usually we omit the symbol $\TrP$ and, for short, write $\tilde f
\tilde g$ for $\tilde f \TrP \tilde g$.  Since $\tilde f$ is
unique (and thus canonical) representative of a class $\tC_f$, we
have $\tilde\Trop[x_1,\dots,x_n] \cong \Trop[x_1,\dots,x_n] / _
\sim $ and therefore get the one-to-one correspondence
$$\tilde\Trop[x_1,\dots,x_n] \To \tF(\Trop^{(n)})$$
between full polynomials and polynomial function.   In the rest of
our exposition we appeal to the reduced domain
$\tilde\Trop[x_1,\dots,x_n]$.

\begin{definition}\label{def:rRed}
A polynomial $\tilde f\in \tTrop[x_1, \dots, x_n]$ is said to be
\textbf{reducible} if $\tilde f = \tilde g \tilde h$ for some
nonconstant  $\tilde g, \tilde h \in \tTrop[x_1, \dots, x_n]$,
otherwise $\tilde f$ called is \textbf{irreducible}. The product $
\tilde f = \tilde q_1 \cdots \tilde q_s$ is called a maximal
\textbf{factorization} of $\tilde f$ \textbf{into irreducibles} if
each of the $\tilde q_i$'s is irreducible. We say that $\tilde g
\in \tTrop[x_1, \dots, x_n]$ \textbf{divides} $\tilde f$  if
$\tilde f =  \tilde q \tilde g$ for some $\tilde q \in \tTrop[x_1,
\dots, x_n]$.
\end{definition}

\noindent  We instantly encounter new difficulties. \pSkip (i) Not
every nonlinear polynomial $\tfF \in \tTrop[x]$ is reducible; for
example one can easily check that $\tfF = x^2 + \uuu{2}x + 3$ is
irreducible. \pSkip
(ii)
 The factorization into irreducibles need not necessarily be
unique; for example $x^2 \TrS \uuu{2} = (x \TrS \uuu{1})^2$ and at
the same time $x^2 \TrS \uuu{2} = (x \TrS 1)(x \TrS \uuu{1})$,
while $x \TrS \uuu{1} \neq x \TrS 1$. \pSkip
(iii) $\ra$ can be a root of a polynomial $f$, but $(x \TrS a)
\nmid f$, for example $1$ is a root of $f = x^2 \TrS x \TrS 2$ but
$(x \TrS 1) \nmid f$.

\begin{proposition} The polynomial $\tilde g$ divides $\tilde f$, i.e. $\tilde g |  \tilde f$, iff the essential part of $\tilde q \tilde g$
is the essential part of $\tilde f$ for some $\tilde q$, which
means $(\tilde f \TrS \tilde q \tilde g)^e$ is ghost.
\end{proposition}
\begin{proof} $\tilde g |  \tilde f$ iff $\tilde f = \tilde g \tilde
q$, for some $\tilde q$, which means $f \sim gq$ for any $f \in
C_{\tilde f}$, $gq \in C_{\tilde g \tilde q}$. Then, by Corollary
\ref{cor:essential} we get $f^e = (gq)^e$.
\end{proof}

\subSecSpc
\subsection{Tropical polynomials in one indeterminate}\label{sec:TropicalFunctionOneVariable}

The use of the reduced domain $\tilde\Trop[x]$ makes the
development of the theory of polynomials in one intermediate quite
close to the classical commutative theory. We start our exposition
with tangible polynomials and then extend the results to whole
$\tilde\Trop[x]$.

\begin{remark}\label{rmk:ceoffCovexity}
Suppose $\al_i$, $\al_j$, $\al_k \in \Real$ are three tangible
coefficients of $f = \bigoplus_i \al_i x^i$ in  $\Trop[x]$, where
$i<j<k$, then $\al_j \in \uCH(\tpF)$ only if
$$\al_j \geq \frac{\al_i \cdot (k-j) + \al_k \cdot (j-i)}{k-i}.$$
(The arithmetic operations here are the classical ones.) This
relation is simply derived form the convexity of $\CH$, and the
fact that $\uCH$ is its upper part.
\end{remark}

\begin{theorem}\label{thm:oneVariableDecomposition}
Any full-tangible polynomial $\tilde f \in \tilde \Trop [x]$ is
factored uniquely into a product of tangible linear polynomials.
\end{theorem}

\begin{proof}
Proof by induction on $n = \deg (\tilde f) $.  Dividing out by
$\al_n$, we may assume that $\tilde f$  monic. The assertion is
obvious for $n=1$. For $n=2$, given $\tfF(x) = x^2 \TrS \al_1 x
\TrS \al_0$, cf. Remark \ref{rem:convexity}, we have:

  $$\tilde f =  \left\{%
\begin{array}{ll}
    (x \TrS \sqrt{\al_0})^2, &  \al_1 \preceq \sqrt{\al_0}; \\ [2mm]
    (x \TrS \frac{\al_1}{\al_0})(x \TrS
  \al_0) , & \al_1 \succ \sqrt{\al_0}. \\
\end{array}%
\right.    $$ (Here, $\sqrt{\al}$ stands for the tropical square
root, which, in the standard meaning, is just $\frac{\al}{2}$ up
to ghost indication.)

Suppose $n > 2$, if $\tilde f  = x^j \tilde g$, for some $j < n$
we are done by the induction assumption. Otherwise
        $$\tilde f = x^n \TrS \al_{n-1} x^{n-1} \TrS \cdots  \TrS  \al_{1} x \TrS
        \al_{0}, $$
        with $\al_0 \neq \tUniS$. Recall that since $f$ is full,
        $\al_i \neq \tUniS$ for all $i=
0,\dots,n$, and each $(i, \al_i)$  appears on $ \uCH(\tilde f )$,
        but $(i, \al_i)$ is not necessarily a vertex.
        We claim that
        $$ \tilde f = (x \TrS \al_{n-1})
        \OP x^{n-1} \TrS \frac{\al_{n-2}}{\al_{n-1}}x^{n-2} \TrS  \cdots  \TrS  \frac{\al_{1}}{\al_{n-1}} x \TrS
        \frac{\al_{0}}{\al_{n-1}} \CP .$$
 This completes the proof by induction.

        To proof
        this we need to show that $\frac{\al_{i-1}}{\al_{n-1}} \prec
        \al_{i}$ for any $i = 1,\dots,n-1$. Recall that $\tilde f$ is monic, that is
        $\al_n =0$. Assume $\frac{\al_{i-1}}{\al_{n-1}}
        \succeq \al_{i} = \frac{\al_i}{\al_n}$, and thus
        $\frac{\al_{n}}{\al_{n-1}} \succeq
        \frac{\al_{i}}{\al_{i-1}}$. If the inequality is equality,
        it contradicts the essentially of $\al_i x^i$ for $\tilde
        f$ (since then, it would be quasi-essential). Otherwise, it contradicts the
        proprieties in which the sequence $M_f$ of the edges' slopes is descending
        (cf. Lemma \ref{rem:desending}).

Conversely, any different products of tangible linear polynomials
clearly produces a different essential complex, and thus the
factorization of a tangible-full polynomial into linear factors is
unique.
\end{proof}

\noindent The above theorem can be implemented in the following
algorithm:
\begin{algorithm}\label{alg}
  \bfem{(Decomposition algorithm) } Let $\tilde
f = \bigoplus_i \al_ix^i$ be a full-tangible polynomial in $\tilde
\Trop[x]$, the algorithm acts recursively: \pSkip
 (i) if $\tfF$ is
not monic set $\tfF^{(1)} = \bigoplus_i (\al_i/\al_n)x^i$
  and apply the algorithm for  $\tfF^{(1)}$, otherwise
  \pSkip (ii) write  $\tfF = (x \TrS \al_{n-1})\tfF^{(1)} = (x \TrS \al_{n-1})
        (x^{n-1} \TrS \frac{\al_{n-2}}{\al_{n-1}}x^{n-2} \TrS  \cdots  \TrS  \frac{\al_{1}}{\al_{n-1}} x \TrS
        \frac{\al_{0}}{\al_{n-1}})$,
\pSkip
   (iii)
  apply the algorithm again to $\tfF^{(1)}$.
\end{algorithm}

The algorithm is applied for full-tangible polynomial, therefore:
\begin{corollary}\label{cor:factorLinear} The factorization of full-tangible polynomials
is unique, in particular, each is factored uniquely into linear
terms.
\end{corollary}

\begin{remark}
Any linear factor of $\tfF$ determines a root of $\tilde f$,
indeed, assume $(x \TrS \ra)$ is a factor of  $\tfF$ then $\tfF =
(x \TrS \ra) \tfG$ and thus $\tfF(\ra) = (\ra \TrS \ra) \tfG(\ra)
\in \etUnit$. The factorization of $\tilde f$ may contain
identical components, is such a case the multiplicity of a root is
defined to be the number of the corresponding (identical)
components in the factorization.
\end{remark}

\begin{example} The algorithm is simulated for $\tfF =  2x^4 \TrS 5x^3 \TrS 5x^2 \TrS 3x \TrS
0$:
\begin{enumerate}
    \item  $\tfF = 2(x^4 \TrS 3x^3 \TrS 3x^2 \TrS 1x \TrS
    (-2))$ \pSkip
    \item  $\quad = 2(x \TrS 3)(x^3 \TrS x^2 \TrS (-2)x \TrS
    (-4))$ \pSkip
    \item $\quad = 2(x \TrS 3)(x^3 \TrS 0)( x^2 \TrS (-2)x \TrS
    (-4))$ \pSkip
    \item $ \quad= 2(x \TrS 3)(x^3 \TrS 0)( x \TrS (-2))(x \TrS
    (-2))$.
\end{enumerate}
Thus, $3$, $0$, and $-2$ are roots of $\tfF$, where $-2$ has
multiplicity 2. Since $\tf$ is full-tangible then the above
factorization to product of linear terms is unique.
\end{example}

Next we look at nontangible full polynomials. Let us call a
polynomial $\tfF = \bigoplus _{i=0}^t \al _i x ^i$
\bfem{semitangible-full} if $\tfF$ is full with $\al_t$ and $\al
_0$ tangible, but $ \al _i$ are ghost for all $0<i<t$. Dividing
out by $\al_t,$ we may assume that any semitangible-full
polynomial is monic.
\begin{observation}\label{nontangfact} Recall that the restriction of the epimorphism $\epiToMaxPlus: \Trop \to
\eReal$ to $\eReal$ is the identity map while for any $\uuu{\ra}
\in \tUnit$ is given by $\epiToMaxPlus(\uuu{\ra}) = \ra$ (see
\Ref{eq:epiTopicalSemiRings}). Suppose $\tfF = \bigoplus _{i=0}^t
\al _i x ^i$ is monic semitangible-full. Then taking $\beta _i=
\epiToMaxPlus( \frac {\al _{t-1}}{\al_i})$, we have
\begin{equation}\label{eq:fact1} \tfF = ( x ^2 \TrS \al _{t-1}x
\TrS \beta _{t-1}) \tfG,
\end{equation} where $\tfG =
x^{t-2} \TrS \bigoplus _{i=1}^{t-3} \beta_i ^\nu x^i \TrS \frac
{\al _0}{\beta _{t-1}}$. Note that this factorization is not
unique; we could factor out any two roots of $\epiToMaxPlus_*(
\tilde f)$ to produce the first factor, just as long as they are
not both maximal or both minimal.

Suppose $\al_t= 0^\nu$, namely $\tilde f$ is not
semitangible-full, and let $\bt = \epiToMaxPlus({\al _{t-1}})$
then
$$\tfF = (x ^\nu \TrS \bt)\bigoplus _{i=0}^{t-1} \frac {\al _i}{\bt} x ^i.$$
\end{observation}

Therefore, whenever the leading terms are ghost we can use
 Observation \ref{nontangfact} to factor out linear factors $(x ^\nu \TrS a)$ until we reach a tangible leading term. But if we do this
twice, we observe for $a \succ b$ that
$$(x^\nu  \TrS a)(x^\nu \TrS b) = 0^\nu x ^2 \TrS a^\nu x \TrS a  b = (x \TrS a)(x^\nu \TrS
b).$$ Thus, we can always make sure that our factorization has at
most one  linear  factor $x ^\nu  \TrS a$ (for $a$ tangible, and
this is the maximal $a$ of those which appear in the linear
factors  $x ^\nu  \TrS a$ ).

 Likewise, when  the constant term is ghost we can factor out some linear
factor  $x \TrS b^\nu$, and arrange for the constant term to be
tangible. Since, in the above notation, $a^\nu \succ b^\nu,$ we
also have
$$(x ^\nu \TrS a)(x \TrS b^\nu) = 0^\nu x ^2 \TrS a x \TrS (a
b)^\nu.$$

 Iterating, we have the following result:
\begin{proposition}\label{fullfact} Every full polynomial is the product
of at most one linear factor of the form $(x^\nu  \TrS a)$, at
most one linear factor of the form $(x \TrS b^\nu),$ and a
semitangible-full polynomial (which can be factored as in
\Ref{eq:fact1}).\end{proposition}

 Putting together Theorem \ref{thm:oneVariableDecomposition} and Observation
\ref{nontangfact}, we see that any irreducible full polynomial
must have no tangible interior vertices, and at most one interior
lattice point (which must be a nontangible vertex), and thus must
be quadratic, of the form $\al_2 x^2 + \al_1^\nu x + \al_0,$ where
$\al_1^\nu x $ is essential. In conjunction with Corollary
\ref{cor:factorLinear} and Proposition \ref{fullfact}, we have
proved the following result:

\begin{theorem}\label{thm:fullfact} Any full polynomial is the unique product of a full tangible
polynomial (which can be factored uniquely into tangible linear
factors), a linear factor $(x^\nu  \TrS a)$, a linear factor $(x
\TrS a^\nu)$, and a semitangible-full
 polynomial, and this factorization is unique.
\end{theorem}
\begin{proof} Just factor at each tangible monomial, and multiply
together the full tangible factors.\end{proof}

\begin{note} Recall that using the (r,u)-decomposition
any polynomial  $\tfF \in \tilde \Trop[x]$ can be written uniquely
as $ \tfF= \tpFr \TrS \tpFu, $ where $\tpFr$ and $\tpFu$ are
tangible polynomials. Using Theorem
\ref{thm:oneVariableDecomposition}, each of these components can
be factored uniquely to a product of linear factors, therefore
$\tfF $ can be written as $$ \tfF(x) = \bigodot_j(x \TrS
\aaa_j)^{i_j} \TrS\bigodot_k(x \TrS \aab_k)^{h_k}
$$
and this decomposition is unique.
\end{note}

\subSecSpc
\subsection{Tropical polynomials in several indeterminates }\label{sec:TropicalFunctionSeveralVariables}

Polynomials in $\tTrop[x_1, \dots, x_n]$ have some special
properties, mainly due to their combinatorial nature. (Recall that
$\bfi = (i_1,\dots,i_n)$ stands for a multi-index and $\bfx  =
(x_1,\dots,x_n)$.)

\begin{proposition}\label{thm:maximumMonomialValue} Let $\tfF= \bigoplus_\bfi \al_\bfi
\bfx^\bfi$ and let $\tfG = \bigoplus_\bfi (\al_\bfi \bfx^\bfi)^k$
for some positive $k \in \Net$. Given $\bfa= (\ra_1,\dots,\ra_n)
\in \Trop^{(n)}$, assume $\tfF(\bfa) = h_\bfi(\bfa)$ for some
monomial $h_\bfi$ of $\tfF$, then $\tfG(\bfa) = (h_\bfi(\bfa))^k$.
\end{proposition}
\begin{proof}
Assume $(\al_\bfi \bfa^\bfi)^k \prec (\al_\bfj \bfa^\bfj)^k$, but
this means $\al_\bfi \bfa^\bfi \prec \al_\bfj \bfa^\bfj$ -- a
contradiction.
\end{proof}

\begin{proposition}\label{thm:maximumMonomialValue.1} For any $\tfF,\tfG \in
\tTrop[x_1,\dots,x_n]$ and any positive  $k \in \Net$, $(\tfF \TrS
\tfG)^k =  \tfF^k \TrS \tfG^k$.
\end{proposition}
\begin{proof} Expand the product $(\tfF \TrS
\tfG)^k$ and observe a mixed component $\tfF^i \tfG^j$, with $i+j
=k$ and $i,j \neq 0$. Pick $\bfa \in \Trop^{(n)}$ and assume
$\tfF(\bfa) \succeq \tfG(\bfa)$, then $\tfF(\bfa)^i  \tfG(\bfa)^j
\preceq \tfF(\bfa)^i \tfF(\bfa)^j = \tfF(\bfa)^k$. On the other
hand, if $\tfF(\bfa) \preceq \tfG(\bfa)$, then $\tfF(\bfa)^i
\tfG(\bfa)^j \preceq \tfG(\bfa)^k$. This means $\tfF^i  \tfG^j$ is
inessential.
\end{proof}

\begin{theorem}\label{thm:powerOfDunctionSeveralVariables} Let $\tfF= \bigoplus_\bfi \al_\bfi
\bfx^\bfi$ and let $\tfG = \bigoplus_\bfi (\al_\bfi \bfx^\bfi)^k$
for some positive $k \in \Net$, then  $\tfF^k = \tfG$.
\end{theorem}
\begin{proof} By the law of polynomial multiplication, it is
clear that as a polynomial $\tfF^k$ has more monomials than $\tfG$
(i.e. all the monomials of $\tfG$ appear also in $\tfF^k$). If
$\tpF$ have a single monomial we are done. Otherwise, pick a
monomial $h_i$ of $\tfF$  and write $\tfF = h_\bfi \TrS \tfF_1$.
Using Proposition \ref{thm:maximumMonomialValue.1}, $\tfF^k =
h_\bfi^k \TrS \tfF_1^k$. Now proceed inductively on $\tfF_1$ to
complete the proof.
\end{proof}

\begin{example}
Let $\tfF(x,y) = x \TrS y$ then, by taking the full closures we
have $$\tfF^2(x,y)= (x \TrS y)^2 = x^2 \TrS \uuu{0}xy \TrS y^2 =
\widetilde{x^2 \TrS y^2}.$$
\end{example}

\secSpc
\section{Tropical Algebraic Sets and Com-sets}\label{sec:TropicalVarieties}
As in the classical theory, using the notion of algebraic sets we
establish the connection between polynomials and tropical
geometry. It turns out that by introducing a new notion of
tropical algebraic com-set the development becomes much easier and
allows the formulation of tropical analogues to classical results,
the tropical Nullstellensatz will be our main example.

Despite our main interest, from the point of view of commutative
algebra, is mainly in the tropical reduced domain (cf. Definition
\ref{def:reduced}), the development in this and in the next
section is being made in the framework of the extended tropical
polynomial semiring  $\Trop[x_1,\dots,x_n]$ that is much wider.

\subSecSpc
\subsection{Tropical algebraic sets}\label{sec:TropVarieties}
\begin{definition}\label{def:TropicalVariety}
The \textbf{tropical algebraic set} of a non empty subset $\tpFF
\subseteq \Trop[x_1,\dots,x_n]$ is defined to be
\begin{equation}\label{eq:tropicalVarieties}
  \afSet(\tpFF) = \{ \bfa \in \Trop^{(n)} \; | \;
  f(\bfa) \in \etUnit,  \; \; \forall f  \in F\}.
\end{equation}
$\afSet(\tpFF)$ is sometimes called  tropical set, for short, and
we call its elements \textbf{roots}, or  \textbf{zeros},  of
$\tpFF$. We say that a subset $\aSet \subset \Trop^{(n)}$ is
\textbf{algebraic}, in the topical sense, if $\aSet =
\afSet(\tpFF)$ for a suitable  $\tpFF \subseteq
\Trop[x_1,\dots,x_n]$.
\end{definition} \noindent
Note that, if $\bfa \in \afSet(\tpFF)$ we necessarily have
$\uuu{\bfa} \in \afSet(\tpFF)$, but the converse claim is not
true.

\begin{remark}\label{rmk:closedAS}In our topology, over closed set, the operations  $\TrS$ and
$\TrP$ are continuous, and the sets $\{ -\infty\}$ and $\etUnit$
are closed. Accordingly, tropical polynomials are continuous as
well, cf. Definition \ref{def:topology}.
\end{remark}

\begin{figure}[!h] \qquad

%
%
%
%
%
%

\begin{minipage}{0.4\textwidth}
\begin{picture}(10,110)(0,0)
\includegraphics[width=\FigWidth in]{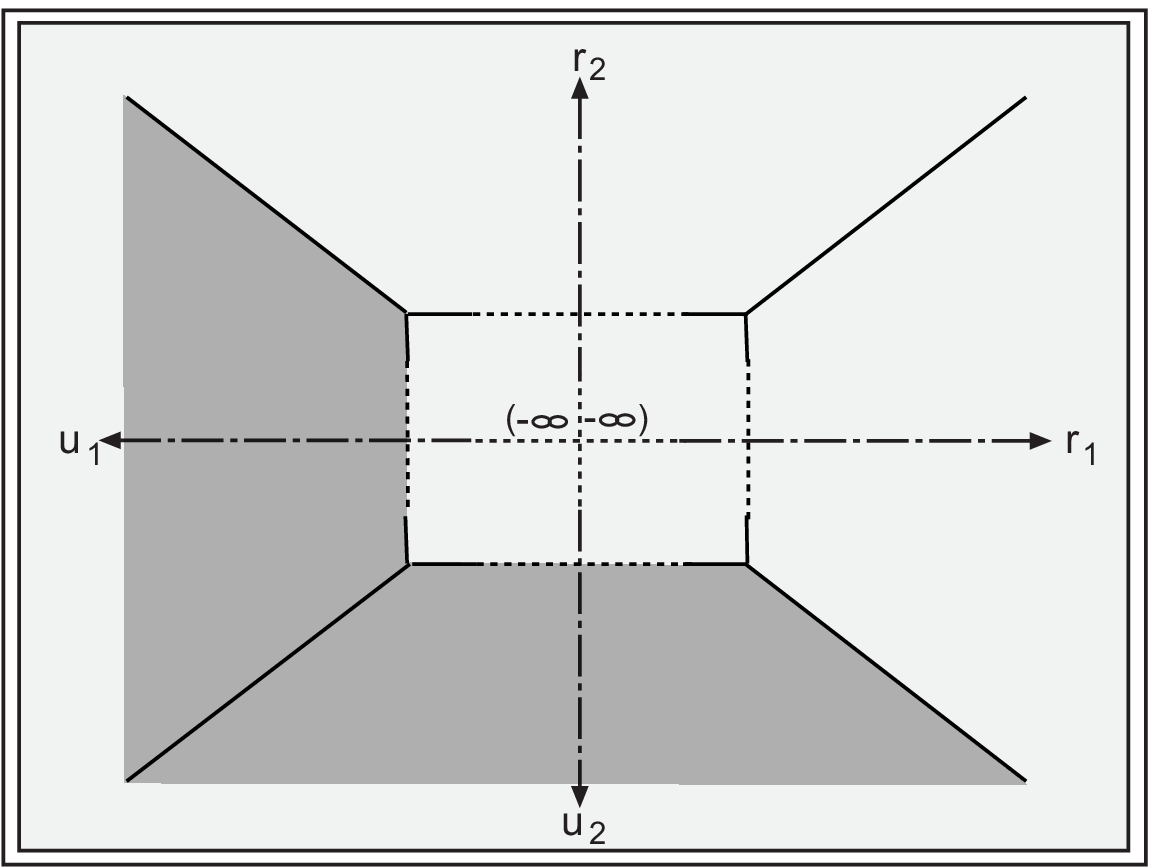}
\end{picture}

\end{minipage}\hfil
\begin{minipage}{0.4\textwidth}
\begin{picture}(10,110)(0,0)
\includegraphics[width=\FigWidth in]{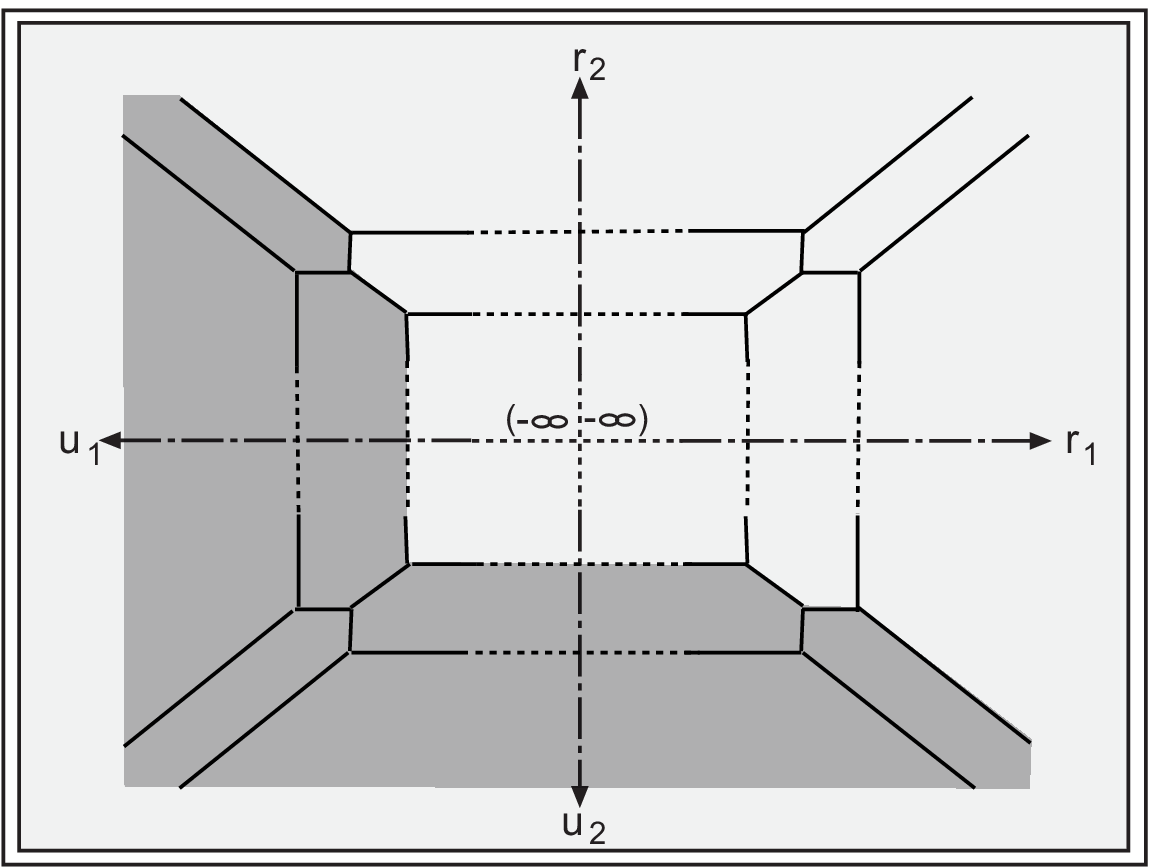}
\end{picture}

\end{minipage}
\caption{\label{fig:tropicalLine} Tropical line and  tropical
conic in $\Trop^{(2)}$.}
\end{figure}

Clearly, for any $f \in \Trop[x]$, $\afSet(f)$ is just the set of
roots of $f$. Analogously, we consider a tropical algebraic set
$\afSet(\tpFF)$ as the set of common solutions of all members of
$\tpFF$. Therefore,  when $\tpF(\bfa) \in \etUnit$, we keep the
familiar terminology and say that $\tpF$ \bfem{vanishes} at
$\bfa$, or equivalently, that $\tpF$ gives value  in $\etUnit$ at
$\bfa$. When $\tpFF$ has a single member, then $\afSet(\tpF)
\subset \Trop^{(n)}$ is called a \bfem{tropical hypersurface}; as
an example see Fig. \ref{fig:tropicalLine}.

\parSpc

\begin{example} Let $\tpF_1 = x_1 \TrS 1$ and  $\tpF_2 = x_2 \TrS 1$
be polynomials in $\Trop[x_1,x_2]$, then
$$ \afSet(\tpF_1) = \{ (1,y) \; | \; y \in \Trop \} \; \cup \;
    \{ (x,y) \; | \; 1 \preceq x \in \etUnit,  y \in \Trop\},$$
while the tropical set of $\tpF_1$ and $\tpF_2$ is the union:
$$ \begin{array}{ll}
  \afSet(\tpF_1,\tpF_2) = & \{ (1,1)\} \; \cup \;
    \{ (1,y) \; | \; 1 \preceq y \in \etUnit\} \; \cup \; \\ [2mm]
   & \{ (x,1) \; | \; 1 \preceq x \in \etUnit\} \; \cup \;
    \{ (x,y) \; | \; 1 \preceq x,y \in \etUnit\}.
  \end{array}
$$
Here, $(1,1)$ is the only common tangible zero.
\end{example}

\begin{lemma}\label{lem:openAlgSet}
Assume $\aSet$ is tropical algebraic set then $\aSet$ is closed
set in the topology of $\Trop^{(n)}$.
\end{lemma}

\begin{proof}
We may assume $\aSet  = \afSet(\tpF)$, for $f = \bigoplus_i f_i$ a
sum of monomials $f_i$'s, is a tropical hypersurface, otherwise
$\aSet = \afSet(\tpFF)$ will be an intersection of closed sets.
Pick a point $\bfa \notin \aSet$ in the complement of $\aSet$,
then we have $f(\bfa) = f_i(\bfa) \in \Real$, for some monomial
$f_i$. Assume first, that all the coordinates of $\bfa$ are
tangible. In the classical sense $f_i$ is smooth and linear, so
there is an open neighborhood $U \subset \Real^{(n)}$ of $\bfa$
such that $f(\bfb) = f_i(\bfb)$ for each $\bfb \in U$. This
implies the complement is open.

If $\bfa$ has a ghost coordinate then $f_i$ is a tangible
constant, since otherwise $\bfa$ would be in $\aSet$, them use the
same argument of the previous paragraph.
\end{proof}

 The next lemma determines the operations on tropical algebraic sets.
\begin{lemma}\label{thm:booleanOpreationInVarities}
Assume  $\aSet',\aSet'' \subseteq \Trop^{(n)}$ are tropical sets,
then so are $\aSet' \cap \aSet''$ and $\aSet' \cup \aSet''$.
\end{lemma}

\begin{proof}
Suppose $\aSet' = \afSet(\tpFF)$ and $\aSet'' = \afSet(\tpGG)$,
where  $\tpFF, \tpGG \subset \Trop[x_1,\dots,x_n]$ are nonempty.
We claim that
$$
  \aSet' \cap \aSet'' =
  \afSet(\tpFF \cup \tpGG) \qquad \text{and} \qquad
  \aSet' \cup \aSet'' =  \afSet(\tpF \tpG \; : \; f \in \tpFF, g
  \in \tpGG).
$$
The left part is by definition; assume $\bfa \in \aSet' \cap
\aSet''$ then $\tpF(\bfa) \in \etUnit$ and $\tpG(\bfa) \in
\etUnit$ for each  $\tpF \in \tpFF$ and  $\tpG \in \tpGG$, which
is the same as all the members of $\tpFF \cup \tpGG$ give values
in $\etUnit$.

For the right part, if $\bfa \in \aSet'$, then all the $\tpF$'s of
$\tpFF$ give values in $\etUnit$ at $\bfa $, which implies that at
$\bfa$ all the products $\tpF\tpG$ also give values in $\etUnit$.
Thus $\aSet' \subset \afSet(\tpF\tpG)$, and $\aSet'' \subset
\afSet(\tpF\tpG)$ follows similarly. This proves the containment
$\aSet' \cup \aSet'' \subset \afSet(\tpF \tpG)$. Conversely,
assume $\bfa \in \aSet(\tpF\tpG)$. If $\bfa \in \aSet'$ we are
done; otherwise $\tpF'(\bfa) \notin \etUnit$ for some $\tpF' \in
\tpFF$, i.e. $\tpF'(\bfa) \in \Real$. But, since at $\bfa$,
$\tpF'\tpG$ gives value in $\etUnit$  for all $\tpG \in \tpGG$,
then $\tpG$ must give value in $\etUnit$ at $\bfa$. This proves
that $\bfa \in \aSet''$, and hence $\afSet(\tpF\tpG) \subset
\aSet' \cup \aSet''$.
\end{proof}

\begin{remark}\label{rem:setofPow}
 From the Lemma and Proposition \ref{thm:maximumMonomialValue} we
 can conclude that $\afSet(f) = \afSet(f^k)$ for each $f \in
 \Trop[x_1,\dots,x_n]$ and any positive $k \in \Net$.
\end{remark}
\begin{remark} \textbf{Tropicalization and tropical
sets:} Based on Kapranov's Theorem \cite{Develin2003,Shustin1278},
the classical tropical hypersurface over $\eMaxPlusAlg$ is the
corner locus (i.e domain of non-smoothness) of a convex piecewise
affine linear function of the form
 \begin{equation}
\nAmb_f  = \max_{\bfi}(\nVal(c_\bfi) + \bfi.\nXx )
\end{equation}
where $\bfi.\nXx$ stands for the standard inner product and the
$c_\bfi$'s are coefficients of a ``superior'' polynomial $f \in
\Fld[\nTz]$ over a non Archimedean field $\Fld$ with a real
valuation $\Val$. Namely, a point $\bfa$ belongs to the corner
locus exactly when two components of $\nAmb_f $ simultaneously
attain the maximum. This is precisely our interpretation of the
tropical addition in view of Definition \ref{def:TropicalVariety}.

In other words, one can consider $\nAmb_f$ as a tangible
polynomial in $\Trop[\nTx]$, then its corner locus with respect to
$\eMaxPlusAlg$, is exactly the restriction of $\afSet(\nAmb_f)$ to
$\Real^{(n)}$, cf. Note \ref{obs:geometryOfGraphs}.
 \end{remark}

\subSecSpc
\subsection{Tropical algebraic com-sets}\label{sec:TropVarieties}
The next object we introduce is  central for our future
development. Given a tropical algebraic set $\aSet \subset
\Trop^{(n)}$ we denote the complement  of $\aSet$ by $\aSet^c$,
that is
$$\aSet^c = \Trop^{(n)} \setminus \aSet \ . $$
Recall that $\aSet$ is a closed set in the topology of $\Trop$, so
for our purpose, connectedness of subsets of $\Trop^{(n)}$ is well
defined.
\begin{definition}\label{def:coVar} Given a tropical algebraic set $\afSet(\tpF) \subset
\Trop^{(n)}$, $\tpF \in \Trop[\nTx] $, the set
\begin{equation}\label{eq:com-set}
\cfSet(\tpF) = \{ \conC_\tpF \; | \;
  \conC_\tpF \text{ is a  connected  component  of } \afSet(\tpF)^c\},
\end{equation}
is defined to be the \textbf{tropical algebraic com-set} (or
tropical com-set, for short) of $\tpF$. A set $\cSet = \{ D_t \; |
\; D_t \subset \Trop^{(n)} \}$ is said to be algebraic com-set if
$\cSet = \cfSet(\tpF)$ for some $\tpF \in \Trop[\nTx] $.
\end{definition}
\noindent Accordingly, any member $\conC_\tpF$  of $\cfSet(\tpF)$
is an open set, cf. Remark \ref{rmk:closedAS}. Since the two
enlarged copies of $\Real$ (i.e. $\eReal$ and $\etUnit$) are glued
along $\tUniS$, the connectivity of components may comprise paths
through $\tUniS$; for instance, the set
$$\conC_\tpF = \{ x \in \Trop \; | \; \uuu{\aaa} \succ x \prec
\aab \}$$ is a proper connected component with $\tUniS \in
\conC_\tpF$.

\begin{example}\label{exmp:coVar}
The tropical algebraic  com-set of  $\tpF = x \TrS \aaa$ is
$$\cfSet(\tpF) = \left\{
\{ x \in \Trop \; | \; \uuu{\aaa} \succ x \prec \aaa \}, \; \{ x
\in \Real \; | \; \aaa \prec x \} \right\}.$$
\end{example}

\begin{remark}
In  view  of  Definition \ref{def:coVar}, over each $\conC_\tpF
\cap \Real^{(n)}$,  $\conC_\tpF \in \cfSet(\tpF)$, $\tpF$ is a
continuous smooth linear function (in the standard meaning), where
for all $\bfa \in \conC_\tpF$, either $\tpF(\bfa) \in \Real$ or
$\tpF(\bfa) \in \tUnit$.
\end{remark}

To emphasize, a tropical com-set is the set of connected
components (each is a set by itself)  of the complement of a
tropical algebraic set. For the forthcoming development, we define
the union
\begin{equation}\label{eq:uCoVar}
\uCoVar(\tpF)= \bigcup_{\conC_\tpF \in \cfSet(\tpF)} \conC_\tpF
\end{equation}
of all the members of $\cfSet(\tpF)$. Therefore, $\uCoVar(\tpF)
\subseteq \Trop^{(n)}$ and  $\afSet(\tpF)^c =  \uCoVar(\tpF).$

\parSpc
\begin{example}\label{exmp:spepcialVar}
Here are some typical cases, assume $\tpF \in \Trop[x]$ then:
\pSkip (i) if $\tpF$ is a tangible constant, i.e. $\tpF \in
\Real$, then $\afSet(\tpF) = \emptyset$ and $\cfSet(\tpF) =
\{\Trop\}$;
\pSkip (ii) if $\tpF$ is a ghost polynomial, $\afSet(\tpF) =
\Trop$ and $\cfSet(\tpF) = \{\emptyset\}$;
\pSkip
 (iii) if $\tpF  = x$
then $\afSet(\tpF) = \etUnit$ and $\cfSet(\tpF) = \{ \Real \}$;
\pSkip (iv) when $\tpF  = \tUniS$ then $\afSet(\tpF) = \Trop$ and
$\cfSet(\tpF) = \{ \emptyset \}$.
\end{example}

We also have the analogous properties to that of tropical
algebraic sets:
\begin{lemma}\label{thm:coVarOfPower} For any $\tpF, \tpG \in
\Trop[\nTx]$: \pSkip (i) $\cfSet(\tpF) = \cfSet(\tpF^k)$;
\pSkip (ii)
 $ \cfSet(\tpF \tpG) = \{ \conC_\tpF \cap  \conC_\tpG \neq \emptyset \; | \;
  \conC_\tpF \in \cfSet(\tpF), \conC_\tpG \in \cfSet(\tpG) \}$; \pSkip
  (iii) for any $\conC_{\tpF \tpG} \in \cfSet(\tpF \tpG)$ there exists
        $\conC_{\tpF} \in \cfSet(\tpF)$ such that $\conC_{\tpF
        \tpG}\subseteq \conC_{\tpF}$.
\end{lemma}
\begin{proof} $(i)$ is obtained directly from the equality $ \afSet(\tpF) =
\afSet(\tpF^k)$ (cf. Remark \ref{rem:setofPow}).  $(ii)$ By
definition:
$$\cfSet(\tpF \tpG) = \{ \conC \; | \;
  \conC \text{  is a  connected component of } \;
  \afSet^c(\tpF\tpG)\},$$
and thus, $\afSet(\tpF\tpG)^c = \afSet(\tpF)^c \cap
  \afSet(\tpG)^c$. Since $\afSet(\tpF\tpG) = \afSet(\tpF) \cup \afSet(\tpG)$, cf. Lemma
  \ref{thm:booleanOpreationInVarities}, then $\afSet(\tpF\tpG)^c$ consists of all
  nonempty intersections of
  connected components from $\afSet(\tpF)^c$ and
  from $\afSet(\tpG)^c$;   $(iii)$ is then obtained
  directly from $(ii)$.
\end{proof}

We generalize Definition \ref{def:coVar} as follows:
\begin{definition}
The tropical algebraic com-set of a nonempty $\tpFF \subseteq
\Trop[x_1,\dots,x_n]$ is defined as
$$\cfSet(\tpFF) = \bigcup_{\tpF \in \tpFF}  \cfSet(\tpF).$$
(This union is not a disjoint union and identical components have
a single instance in $\cfSet(\tpFF)$.) We say that $\cSet
\subseteq \Trop^{(n)}$ is tropical algebraic com-set if $\cSet =
\cfSet(\tpFF)$ for a suitable $\tpFF \subseteq
\Trop[x_1,\dots,x_n]$.

\end{definition}

\begin{definition} Given tropical algebraic com-sets $\cSet',
\cSet'' \subset \Trop^{(n)} $ we define the intersection $\capS$
to be
\begin{equation}\label{eq:specialIntersection}
  \cSet' \capS \; \cSet'' =
  \{ \conC' \cap  \conC'' \neq \emptyset \; | \;
  \conC' \in \cSet', \conC'' \in \cSet'' \}.
\end{equation}
The inclusion  $\subS$ is defined by the rule:
\begin{equation}\label{eq:specicalContainment}
  \cSet' \subS \cSet'' \; \IFF
 \text{ for each } \conC' \in \cSet' \text{ there exists } \conC'' \in
  \cSet'' \text{ such that } \conC'  \subseteq \conC''.
\end{equation}

\end{definition}

\begin{lemma}\label{thm:booleanOpreationInCoVarities}
Assume  $\cSet',\cSet'' \subset \Trop^{(n)}$ are tropical
algebraic com-sets, then so are $\cSet' \cup \; \cSet''$ and
$\cSet' \capS \; \cSet''$.
\end{lemma}
\begin{proof}
Suppose $\cSet' = \cfSet(\tpFF)$ and $\cSet'' = \cfSet(\tpGG)$,
$\tpFF,\tpGG \subset \Trop[x_1,\dots,x_n]$, are not empties, we
claim that
$$
  \cSet' \cup \cSet'' =  \cSet(\tpFF \cup
  \tpGG), \qquad and \qquad
  \cSet' \capS \cSet'' =  \cfSet(\tpF\tpG\; : \; f\in
  \tpFF,
  g\in \tpGG).
$$
Indeed, the left equality is by definition while the right is the
generalization of Lemma \ref{thm:coVarOfPower} in terms of
Equation \eqref{eq:specialIntersection}.
\end{proof}

\secSpc
\section{Tropical Ideals}\label{sec:TropicalIdeals}
Ideals are  main structure in the classical theory; we develop
this notion in the tropical sense. As will be seen, the tropical
ideal is an analogous of the classical one. Later, we will study
the main properties of tropical ideals and realize how they relate
to tropical sets and com-sets.

\subSecSpc
\subsection{Definition and properties}\label{sec:TropVarieties}
\begin{definition}\label{def:tropIdeal} A subset $\tIdA \subset
\Trop[\nTx]$ is a \textbf{tropical ideal} of polynomials if it
satisfies: \pSkip (i) $\tUniS \in \tIdA$; \pSkip (ii)
  if $\tpF,\tpG \in \tIdA$, then $\tpF \TrS \tpG \in \tIdA$; \pSkip (iii)
   if $\tpF \in \tIdA$,  and $\tpH \in \Trop[\nTx]$, then
  $\tpH\tpF \in \tIdA$. \pSkip
A ideal is called \textbf{tangible ideal} if all of its elements
are tangible and is called \textbf{ghost ideal} when all of its
elements are ghost.
\end{definition}
\noindent As an example, one can easily verify that
$\etUnit[\nTx]$ is a proper tropical ideal of $\Trop[\nTx]$. (Note
that we may have ideal which are neither, tangible ideal nor ghost
ideal.)

An immediate  conclusion is:
\begin{corollary}\label{thm:singleMaximalIdeal}
There exists only a single proper maximal ideal $\mtId \subset
\Trop[\nTx]$.
\end{corollary}
\begin{proof}
We identify the maximal ideal as $\mtId = \Trop[\nTx] \setminus
\Real$, that is the set  of all polynomials in $n$ indeterminate
$x_1,\dots,x_n$ except constant tangible polynomials. Assume that
$\mtId$ can be enlarged further, say by $\aaa \in \Real$. Now, if
$\aaa \in \mtId$ then $(-\aaa) \in \mtId$, and hence $0 \in
\mtId$, which is the multiplicative unit of $\Trop[\nTx]$. But
then, for any $\tpF \in \Trop[\nTx]$ we have $0\tpF = \tpF$ which
means that $\tpF \in \mtId$, thus $\mtId$ is no more a proper
ideal. Clearly, for any other proper ideal $\tIdA \subset
\Trop[\nTx]$ we have $\tIdA \subseteq \mtId$ since otherwise
$\tIdA$ must contain a constant tangible polynomial and by the
previous argument it would not be proper.
\end{proof}

The operations between ideals and a polynomial $\tpF \in
\Trop[x_1,\dots, x_n]$ are defined in terms of elements:
$$
  \tpF \TrS \tIdA  =  \{\tpF \TrS \tpG \; | \; \tpG \in \tIdA  \}
  \qquad and \qquad
  \tpF \TrP \tIdA  =  \{\tpF \tpG \; | \; \tpG \in \tIdA  \}.
$$
Clearly, from the latter operation we have $\tpF \TrP \tIdA
\subset \tIdA$ for any $\tpF \in \Trop[x_1,\dots, x_n]$. The first
natural construction of an ideal is the ideal generated by a
finite number of polynomials.

\begin{definition}\label{def:tropPolyIdeal}
Let $\tpF _1,\dots,\tpF _s$ be a collection of polynomials in
$\Trop[\nTx]$, then we set
$$\idl{{\tpF _1,\dots,\tpF_s}} = \left\{ \bigoplus_i \tpH _i \tpF _i \; | \;
  \tpH _1,\dots,\tpH _s \in \Trop[\nTx]\right\}$$
to be the \textbf{ideal generated by} $\tpF_1,\dots,\tpF_s$. When
 $s=1$ the ideal is
called \textbf{principal ideal}. Given an ideal $\tIdA \subset
\Trop[\nTx]$ we say that  $\tIdA$ is \textbf{finitely generated}
if there exist $\tpF _1,\dots,\tpF _s \in \Trop[\nTx]$ such that
$\tIdA = \langle \tpF _1,\dots,\tpF_s \rangle$, or equivalently,
we say that $\tpF _1,\dots,\tpF _s$ are the \textbf{tropical
generating set} of $\tIdA$.
\end{definition} \noindent
As in the classical case, a tropical ideal may have many different
generating sets.

\begin{claim}
 The set $\idl{{\tpF _1,\dots,\tpF_s}}$ is indeed a
 tropical ideal.
\end{claim}
\begin{proof}
 $\tUniS \in \idl{{\tpF
_1,\dots,\tpF_s}}$ since $\bigoplus_i \tUniS \tpF _i = \tUniS$.
Suppose $\tpF = \bigoplus_i \tpP_i\tpF_i$, $\tpG = \bigoplus_i
\tpQ_i\tpF_i$ and let $\tpH \in \Trop[\nTx]$. Then, using the
polynomial rules, the equations
$$\tpF \TrS \tpG = \bigoplus_i (\tpP_i \TrS \tpQ_i)\tpF_i, \qquad
  \tpH \tpF = \bigoplus_i (\tpH \tpP_i)\tpF_i$$
complete the proof.
\end{proof}

Given an ideal $\tIdA \subset \Trop[\nTx]$, as has been done
previously for subsets of $\Trop[\nTx]$, we define the tropical
algebraic set of $\tIdA$ to be
\begin{equation}\label{eq:varietyOfIdeal}
 \afSet(\tIdA)= \{ \bfa \in \Trop^{(n)} \; | \; \tpF(\bfa) \in \etUnit, \forall \tpF \in
 \tIdA\}.
\end{equation}

\begin{proposition} For any tropical ideals $\tIdA
\subseteq \tIdB$ we have the inverse inclusion $\afSet(\tIdB)
\subseteq \afSet(\tIdA)$.
\end{proposition}
\noindent The proof is technically straightforward, so we omit the
proofs' details.

Earlier, we have shown how tropical sets are obtained from ideals,
but we also have the converse direction in which tropical
algebraic sets give rise to ideals.

\begin{definition}
The ideal of a tropical algebraic set $\aSet \subseteq
\Trop^{(n)}$ is defined to be
$$
\tId(\aSet) = \{ \tpF \in \Trop[\nTx] \; | \; \tpF(\bfa) \in
\etUnit, \; \forall \bfa \in\aSet \}.
$$
\end{definition}
\noindent The crucial observation is that $\tId(\aSet)$ is indeed
a tropical ideal.

\begin{lemma}
Let $\aSet \subset \Trop^{(n)}$ be a tropical algebraic set, then
$\tId(\aSet) \subset \Trop[\nTx]$ is a tropical ideal.
\end{lemma}
\begin{proof}  $\tUniS
\in \tId(\aSet)$ by definition. Assume  $\tpF,\tpG \in
\tId(\aSet)$,  $\tpH \in \Trop[\nTx]$, and  $\tUniS \neq \bfa \in
\aSet$; then
$$
  (\tpF \TrS \tpG)(\bfa) = \tpF(\bfa) \TrS \tpG(\bfa) =
   \uuu{x} \TrS \uuu{y}   \in  \etUnit, \qquad
   (\tpH \tpF)(\bfa) =  \tpH(\bfa) \tpF(\bfa) =  (\tpH(\bfa))\uuu{x}  \in
   \etUnit,
$$ where $\tpF(\bfa) = \uuu{x}$ and  $\tpG(\bfa) = \uuu{y}$,
and it follows that $\tId(\aSet)$ is an ideal.
\end{proof}

\begin{lemma}
Let $\tpF_1,\dots,\tpF_s \in \Trop[\nTx]$, then $\idl{{ \tpF
_1,\dots,\tpF_s }} \subset \tId(\afSet( \tpF _1,\dots,\tpF_s)).$
\end{lemma}
\begin{proof}
For $\tpF \in \idl{{\tpF_1,\dots,\tpF_s}}$ we have $\tpF =
\bigoplus_i, \tpH_i\tpF_i$ where the $\tpH_i$'s are polynomials in
$\Trop[\nTx]$. Since all $\tpF_1,\dots,\tpF_s$ give values in
$\etUnit$ on $\afSet( \tpF _1,\dots,\tpF_s)$, so does $f=
\bigoplus_i, \tpH_i\tpF_i$, which proves that $\tpF \in
\tId(\afSet( \tpF _1,\dots,\tpF_s))$.
\end{proof}

\begin{proposition}
Let $\aSet'$ and $\aSet''$ be tropical algebraic sets then, \pSkip
(i) $\aSet' \subset \aSet''$ if and only if $\tId(\aSet') \supset
  \tId(\aSet'')$; \pSkip (ii)
   $\aSet' = \aSet''$ if and only if $\tId(\aSet') =
  \tId(\aSet'')$.
\end{proposition}
\begin{proof}
$(i)$ Suppose $\aSet' \subset \aSet''$, then any polynomial that
gives value in $\etUnit$ on $\aSet''$ must also give value in
$\etUnit$ on $\aSet'$. This proves $\tId(\aSet') \supset
\tId(\aSet'')$. Assume that $\tId(\aSet') \supset \tId(\aSet'')$,
we know that $\aSet''$ is the tropical algebraic set defined by a
 set $\tpGG \subset \Trop[\nTx]$, and it follows
that $\tpG \in \tId(\aSet'') \subset \tId(\aSet')$ for any $\tpG
\in \tpGG$. Hence, the $\tpG$'s give values in $\etUnit$ on
$\aSet'$. Since $\aSet''$ consists of all common solutions  of the
$\tpG$'s, it follows that $\aSet' \subset \aSet''$. $(ii)$ is an
immediate consequence of $(i)$.
\end{proof}

\parSpc
Earlier we showed that a tropical set determines a tropical ideal,
next we will show that the same is also valid for com-sets. But
first, let's specify the tropical com-set of an ideal. Given an
ideal $\tIdA \subset \Trop[x_1,\dots,x_n]$, its tropical algebraic
com-set is defined to be the set of connected components
$$\cfSet(\tIdA) = \{\conC_\tpF \; | \; \tpF \in \tIdA\} = \bigcup_{f \in \tIdA } \cfSet(\tpF).$$
Defining $\uCoVar(\tIdA) = \bigcup_{\conC_\tIdA \in \cfSet(\tIdA)}
\conC_\tIdA $, we have  $\afSet(\tIdA) = \Trop^{(n)} \setminus
\uCoVar(\tIdA)$; namely,  $\afSet(\tIdA) = \bigcap_{\tpF \in
\tIdA}(\Trop^{(n)} \setminus \uCoVar(\tpF))$.

We also have the converse direction in which tropical algebraic
com-sets give rise to ideals.

\begin{definition}\label{def:coVar}
Let $\cSet$ be a tropical algebraic com-set, the tropical ideal
$\tId(\cSet)$ is defined as
$$\tId(\cSet) = \{ \tpF \in \Trop[\nTx]  \; | \;
  \forall \conC_\tpF \in \cfSet(\tpF), \exists \conC_o \in \cSet \; s.t.
  \;
  \conC_\tpF \subseteq \conC_o  \}. $$
\end{definition}

\begin{proposition}
$\tId(\cSet)$ is indeed a tropical ideal.
\end{proposition}
\begin{proof}Whether $\cSet$ contains a nonempty set or not, i.e. $\cSet = \{\emptyset\}$, $\emptyset \subseteq \conC$ for
  any $\conC \in \cSet$. Since $\cfSet(\tUniS) = \emptyset$, by Example \ref{exmp:spepcialVar}, we
  have $\tUniS \in \tId(\cSet)$.

   Given $\tpF, \tpG \in \tId(\cSet)$, we need to show that
  $\tpF \TrS \tpG \in \tId(\cSet)$. By the way contradiction, assume
  $\tpF \TrS \tpG \notin \tId(\cSet)$; this means  there exists
  $\conC_o \in \cfSet(\tpF \TrS \tpG)$ that is not
  contained in any member of $\cSet$.  Clearly, $\conC_o \cap
  \conC_{\tpF} \neq \emptyset$ or $\conC_o \cap
  \conC_{\tpG} \neq \emptyset$ for some $\conC_\tpF \in \cfSet(\tpF)$ or
  $\conC_\tpG \in \cfSet(\tpG)$, otherwise $(\tpF \TrS \tpG)(\bfa) \in
  \etUnit$ for all $\bfa \in \conC_o $.

Denote the closure of $\conC_{\tpF}$ by $\overline{\conC}_{\tpF}$
and let $\bound \overline{\conC}_{\tpF}$ be the
  boundary of $\overline{\conC}_{\tpF}$.
  Suppose $\conC_o \cap \conC_{\tpF} \neq \emptyset$, then there
  is some
   $\bfa \in \conC_o \cap
  \bound\overline{\conC}_{\tpF}$, and in particular  $\bfa \in \afSet(\tpF)$. But
  then, there is  $\conC_{\tpG} \in \cfSet(\tpG)$, $\conC_o \cap
  \conC_{\tpG} \neq \emptyset$,  such that  $ \bfa \in \conC_o \cap
  \conC_{\tpG}$. Now, since $\conC_o \nsubseteq \conC_{\tpG}$,
  there exits $\bfb \in \conC_o \cap
  \bound\overline{\conC}_{\tpG}$, and thus  $\bfb \in \afSet(\tpG)$.
  Moreover, the intersection  $\bound \overline{\conC}_{\tpF} \cap \bound \overline{\conC}_{\tpG} \neq \emptyset$
  is contained  in  $\afSet(\tpF) \cap \afSet(\tpG)$, so we
  necessarily have
  $\bound \overline{\conC}_{\tpF} \cap \bound \overline{\conC}_{\tpG} \cap \conC_o
  \neq \emptyset$. Therefore, there is $\bfc \in \conC_o$ on which both
  $\tpF$ and $\tpG$ give values in $\etUnit$  -- a contradiction.
    The last condition in which if  $\tpF\in \tId(\cfSet)$ and $\tpH \in
  \Trop[\nTx]$, then $\tpF \tpH \in \tId(\cSet)$ is derived
  immediately as a result of Lemma \ref{thm:coVarOfPower}.
\end{proof}

\begin{lemma}
Let $\tpF_1,\dots,\tpF_s \in \Trop[\nTx]$, then $\idl{{ \tpF
_1,\dots,\tpF_s }} \subset \tId(\cfSet( \tpF _1,\dots,\tpF_s)).$
\end{lemma}
\begin{proof}
For $\tpF \in \idl{{\tpF_1,\dots,\tpF_s}}$, we have $\tpF =
\bigoplus_i, \tpH_i\tpF_i$ with $\tpH_i \in \Trop[\nTx]$.  $\tpF$
is smooth and linear (in the usual sense) on  any $\conC_\tpF \in
\cfSet(\tpF)$ and is equal to $\tpH_j\tpF_j$ for some $1 \leq j
\leq s$. Hence, $\conC_\tpF \in \cfSet(\tpH_j\tpF_j)$ and there is
$\conC_{\tpF_j} \in \cfSet(\tpF_j)$ such that $\conC_\tpF
\subseteq \conC_{\tpF_j}$, cf. Lemma \ref{thm:coVarOfPower}.
\end{proof}

\begin{proposition}
Let $\cSet'$ and $\cSet''$ be tropical algebraic com-sets, then
\pSkip (i) $\cSet' \subS \cSet''$ if and only if $\tId(\cSet')
\subseteq
  \tId(\cSet'')$; \pSkip (ii)
   $\cSet' = \cSet''$ if and only if $\tId(\cSet') =
  \tId(\cSet'')$.
\end{proposition}
\begin{proof}
Suppose  $\cSet' \subS \cSet''$ and  $\tpF \in \tId(\cSet')$, then
for each  $\conC_\tpF \in \cfSet(\tpF)$ there exists $\conC' \in
\cSet'$ such that $\conC_\tpF \subseteq \conC'$. Since $\cSet'
\subS \cSet''$ then there is $\conC'' \in \cSet''$ such that
$\conC' \subseteq \conC''$, and in particular $\conC_\tpF
\subseteq \conC''$, hence $\tpF \in \tId(\cSet'')$.
Conversely, assume \Ref{eq:specicalContainment} is not satisfied,
then there exists $\conC'_o \in \cSet'$ such that $\conC'_o
\nsubseteq \conC''$ for any $\conC'' \in \cSet''$. We know that
$\cSet' = \cfSet(\tpFF) $ for some  $\tpFF \subset \Trop[\nTx]$,
thus $\conC'_o \in \cfSet(\tpF)$ for some $ \tpF \in \tpFF $. In
particularly, $\tpF \in \tId(\cSet')$. But, since $\conC'_o$ is
not contained in any member of $\cSet''$, then $\tpF \notin
\tId(\cSet'')$.
\end{proof}

\subsection{Radical ideals}\label{sec:TropVarieties}
We turn to deal with special types of  ideals.
\begin{definition}\label{def:radicalIdeal}
 The \textbf{radical}  $\rtIdA$ of an ideal $\tIdA \subset
\Trop[\nTx]$ is defined to be the set of all $\tpF \in
\Trop[\nTx]$ for which  $\tpF^k \in \tIdA$ for some positive $k
\in \Net$. An ideal $\tIdA$ is called a \textbf{radical ideal} if
$\rtIdA = \tIdA$. An ideal $\ptId \subset \Trop[\nTx]$ is said to
be \textbf{prime} ideal if when $\tpF \tpG \in \ptId$, then either
$\tpF \in \ptId$ or $\tpG \in \ptId$.
\end{definition}

Any ideal $\tIdA$ is contained in some prime ideal $\ptId$. We can
simply complete it to prime ideal: whenever an element $\tpH
=(\tpF\tpG) \in \tIdA$ and both $\tpF$ and $\tpG$ are not in
$\tIdA$, add one of them (including its multiples) to $\tIdA$. By
this construction, $\tIdA$ is completed to be a prime ideal
$\ptId$. We can conclude that:
\begin{corollary}\label{thm:nilRadicalAndPrimal} Every topical prime ideal is a tropical radical ideal.
\end{corollary}

The next two propositions are immediate.
\begin{proposition} The radical of a tropical ideal $\tIdA$ is again a tropical ideal.
\end{proposition}

\begin{proof} Suppose $\tpF, \tpG \in \rtIdA$, thus $\tpF^k \in \tIdA$ and $\tpG^m \in \tIdA$ for some
positive integers $k,m$. Then
$$
(\tpF \TrS \tpG)^{k+m} = \bigoplus^{k+m}_{i=0} h_i \tpF^i
\tpG^{k+m-i},
$$
where  $h_i \in \Trop[\nTx]$. In each term either $i \geq k$ or $k
+m - i \geq n$. In the first case, $\tpF^i \in \tIdA$, and in the
second case, $\tpG^{k+m -i} \in \tIdA$. Since $\Trop[\nTx]$ is
commutative and $\tIdA$ is an ideal, the sum of these terms is
again in $\tIdA$, and hence $\tpF \TrS \tpG \in \tIdA$. To see
that $\tIdA$ is closed under multiplication by elements
$\Trop[\nTx]$; let $\tpH \in  \Trop[\nTx]$, then $(\tpH \tpF)^m =
\tpH^m \tpF^m \in \tIdA$, so $\tpH \tpF \in \rtIdA$.
\end{proof}

\begin{proposition} The radical of $\rtIdA$ is equal to $\rtIdA$.
\end{proposition}
\begin{proof} Clearly, $\rtIdA$ is contained in the radical of $\rtIdA$. To see the reverse
inclusion, assume $\tpF \in \rtIdA$, then $\tpF^k \in \rtIdA$ for
a positive $k \in \Net$, which means that $(\tpF^k)^m \in \rtIdA$
for some positive $m \in \Net$. Since $\tpF^{km} \in \rtIdA$, we
see that $\tpF \in \rtIdA$.
\end{proof}


\begin{definition}
A polynomial $\tpF \in \Trop[x_1,\dots,x_n]$ is called
\textbf{ghost-potent} if $\tpF^k \in \tUnit[x_1,\dots,x_n]$ for
some positive $k \in \Net$.  A \textbf{ghost-radical} of a ghost
ideal $\tIdA \subset \Trop[\nTx]$ is defined to be the set of all
$\tpF \in \Trop[\nTx]$ for which  $\tpF^k \in \tIdA$ for some
positive $k \in \Net$. The ghost-radical of $\tUnit[\nTx]$,
denoted $\rad(\tUnit)$, is the set of all ghost-potent elements in
$\Trop[\nTx]$.
\end{definition}
\noindent (We may extend this definition by joining  $\tUniS$ to
the ghost ideal.)

Clearly, any ghost element is ghost-potent. Restricting ourselves
to the reduced domain $\tTrop[x_1,\dots,x_n]$ we have the
following:
\begin{proposition}\label{thm:nilRadical}
The  ghost-radical $\rad(\tUnit)$ is unique and is equal to
$\ttUnit[\nTx]$.
\end{proposition}
\begin{proof} Assume, $\tpF \notin \ttUnit[x_1,\dots,x_n]$ and $\tpF^k \in
\ttUnit[x_1,\dots,x_n]$ for some positive $k \in \Net$. Since
$\tpF \notin \ttUnit[x_1,\dots,x_n]$, it has at least one
essential tangible monomial   $\al_\bfi \bfx^\bfi$. But then
$\tpF^k$ has also at least one essential tangible monomial,
specifically $(\al_\bfi \bfx^\bfi)^k$ (cf. Theorem
\ref{thm:powerOfDunctionSeveralVariables}) -- a contradiction
($\tpF^k \notin \ttUnit[x_1,\dots,x_n]$).
\end{proof}
\noindent This statement is true only for the reduced domain,
because we ``ignore'' inessential monomials, and is not true for
the non-reduced polynomial semiring $\Trop[x_1,\dots,x_n]$, take
for instance $\uuu{0}x^2 \TrS x \TrS \uuu{0}$ is not ghost while
$(\uuu{0}x^2 \TrS x \TrS \uuu{0})^k$, for each $k \in \Net$,  is
ghost.

\begin{theorem}\label{thm:radicalAsPrimeIntersection} Let $\tIdA \subset \Trop[\nTx]$ be an
ideal, and let $P$ be the set of all prime ideals $\ptId \supseteq
\tIdA$, then
$$ \rtIdA = \bigcap_{\ptId \in P} \ptId.$$
In particular, $\rad(\tUnit)$ is the intersection of all prime
ideals in $\Trop[\nTx]$ that contain $\tUnit[\nTx]$.
\end{theorem}

\begin{proof} Denoting $\tP_\cap =
\bigcap_{\ptId \in P} \ptId$, we show $\rtIdA = \tP_\cap$ by cross
inclusion.

$(\subseteq)$  Let $\tpF \in  \rtIdA$,  that is $\tpF^k \in \tIdA$
for some positive integer $k \in \Net$, and take $k$ to be the
least $k$ for which this is true. Let $ \ptId \subset \Trop[\nTx]$
be a prime ideal containing $\rtIdA$,  then $\tpF^k \in \ptId$.
Write $\tpF^k = \tpF \tpF^{k-1}$, since $\ptId$ is prime then
either $\tpF$ in $\ptId$ or $\tpF^{k-1}$ in $\ptId$. If $\tpF \in
\ptId$ we are done, otherwise $\tpF \notin \ptId$ and thus
$\tpF^{k-1} \in \ptId$. But this contradicts the assumption that
$k$ is minimal. To see this, just repeat the decomposition
inductively to obtain $\tpF^2 \in \ptId$, namely $\tpF \in \ptId$
-- a contradiction. Thus, $\tpF$ is contained in every prime ideal
$ \ptId \supseteq \rtIdA$.

$(\supseteq)$ We show that if $\tpF \notin \rtIdA$, then there
exists $\ptId \subset \Trop[\nTx]$ such that $\tpF \notin \ptId$
and hence $\tpF \notin \tP_\cap$. This will be done by
constructing a prime ideal that does not contain $\tpF$. Let $\tpF
\in \Trop[\nTx]$ such that $\tpF \notin \rtIdA$, since $\tUniS \in
\rtIdA$, then $\tpF \neq \tUniS$. Let $\tS$ be the family of
ideals of $\Trop[\nTx]$ that do not contain any power of $\tpF$
and do contain $\tIdA$. This family $\tS$ is not empty because
$\tIdA \in \tS$. Also, we see that chains of ideals in $\tS$ have
upper bounds because if $\tpF^k$ is not in any ideal of a given
chain, then it is also not in the union of the ideals in that
chain. So, we can now apply Zorn's Lemma to see that there is some
maximal element $\ptId_{(max)}$ of $\tS$. Since $\ptId_{(max)}$ is
in $\tS$, $\ptId_{(max)}$ does not contain $\tpF$.

We will now show that $\ptId_{(max)}$ is prime. By way of
contradiction, assume  $g,h \in \Trop[\nTx]$ are not in
$\ptId_{(max)}$ but such that $gh \in \ptId_{(max)}$. Since
$\ptId_{(max)}$ is a maximal element of $\tS$, we see that for
some positive integers $k$ and $m$, $\tpF^k \in (g) \TrS
\ptId_{(max)}$ and $\tpF^m \in (h) \TrS \ptId_{(max)}$. But then
$\tpF^{k+m} \in (gh) \TrS \ptId_{(max)} = \ptId_{(max)}$,
contradicting the fact that $\ptId_{(max)} \in \tS$. Thus
$\ptId_{(max)}$ is indeed a prime ideal, and so $\tpF \notin
\tP_\cap$.
\end{proof}

\section{An Algebraic Tropical Nullstellensatz}\label{sec:TropicalNullstellensatz}

All our previous development leads to the foundation  of an
algebraic tropical nullstellensatz --  the weak version and the
strong version; the weak version is stated in terms of both,
tropical algebraic sets and tropical algebraic com-sets, while the
strong version is phrased only in terms of tropical algebraic
com-sets. The latter is an algebraic rephrasing, enabled due to
our semiring structure, of the tropical nullstellensatz that
appeared in \cite{shustinNull} and part of our development is
based on this theorem.

\subsection{Weak Nullstellensatz}

\begin{theorem}\label{thm:commonZero} Let $\tpF_1,\dots,\tpF_s \in \Trop[\nTx]$
 be nonconstant polynomials, then $\afSet(\tpF_1,\dots,\tpF_s)
\neq \emptyset$.
\end{theorem}
\noindent In fact we can also allow constant ghost polynomials,
but the tropical algebraic set of these polynomial is
$\Trop^{(n)}$.
\begin{proof} Suppose $n=1$. For each $i=1,\dots,s$, assume $h_i = \al_i x^i$ is the least
significant monomial of $\tpF_i$; dividing out by
$\epiToMaxPlus(\al_i) x^i$, then $\tpF_i$ has the form
$$\tpF_i(x) = \underbrace{\al'_n x^n \TrS \cdots \TrS
\al'_1 x}_{=\tpG_i(x)} \TrS \bt, \qquad \bt \in \{0, \uuu{0} \}.
$$
We may assume $\tpG_i$ is nonconstant, since otherwise $\tpF_i$
has a single nonconstant monomial, which means that
$\afSet(\tpF_i) = \etUnit$. According to Lemma \ref{lem:cl1}, for
each $i$, there is $r_i \in \Trop$ for which $\tpG_i(r_i) =
\uuu{0}$. Take $r$ to be the ghost of the maximal $r_i$'s, then,
for each $i$, $\tpG_i(r) = \uuu{a}_i \succeq 0$ and  $\tpG_i(r)
\TrS 0 \in \tUnit$. The generalization to $n > 1$ is obvious, just
pick $a \in \Real$, fix $x_2 = \cdots = x_n = a$, and apply the
above argument for $x_1$.
\end{proof}

\begin{corollary}\label{thm:emptyVar}
Let $\tpF_1,\dots,\tpF_s \in \Trop[\nTx]$, then
$\afSet(\tpF_1,\dots,\tpF_s) = \emptyset$ if and only if one of
the $\tpF_i$'s is a constant tangible, i.e. $\tpF_i = \aac \in
\Real$.
\end{corollary}
The corollary is  derived directly from Theorem
\ref{thm:commonZero}.

\begin{theorem}\label{thm:WeakTropicalNullstellensatz}
\textbf{(Weak  Nullstellensatz)} Let $\tIdA \subset \Trop[\nTx]$
be a proper finitely generated ideal, then $\afSet(\tIdA) \neq
\emptyset$. Equivalently, if $\afSet(\tIdA) = \emptyset$, then
$\tIdA = \Trop[\nTx]$.
\end{theorem}

\begin{proof}
Assume $\afSet(\tIdA) = \emptyset$, by Corollary
\ref{thm:emptyVar} there exists a constant tangible polynomial
$\tpF \in \tIdA$, i.e. $f = a \in \Real$.  Then,  $a^{-1} = 0/a
\in \Trop$ and thus $0 \in \tIdA$, which means $0 \tpG = \tpG \in
\tIdA$ for each $\tpG \in \Trop[\nTx]$. This shows that $\tIdA =
\Trop[\nTx]$.
\end{proof}

\begin{corollary}
 Let $\tIdA \subset \Trop[\nTx]$ be a
tropical ideal, then $\Trop^{(n)} \in \cfSet(\tIdA)$ if and only
if $\tIdA = \Trop[\nTx]$.
\end{corollary}
\begin{proof} Immediate, by Theorem
\ref{thm:WeakTropicalNullstellensatz} and the relation:
$\afSet(\tIdA) = \emptyset$ if and only if $\Trop^{(n)} \in
\cfSet(\tIdA)$.
\end{proof}

%

\subsection{Strong Nullstellensatz}
The use of the reduced tropical domain $\tTrop[x_1,\cdot,x_n]$
allows us an easy algebraic formulation of geometric ideas, which
lead to the tropical Nullstellensatz.

\begin{remark}\label{rem:con} Let $f = \bigoplus_i f_i$ and assume $D_f \in \cfSet(\tpF)$.
Then,  $\tpF|_D = \tpF_i|_{D_f}$  for some monomial $f_i= \al_i
x_1^{i_1} \cdots x_n^{i_n}$. Suppose $i_t =0$, for some  $t =
1,\dots,n$, then if the $t$'th coordinate of a point $\bfa \in
D_f$ has a tangible value $a_t$ then, by the connectedness of
$D_f$ the point $\bfa'$, obtained by replacing the coordinate
$a_t$ by $a^\nu_t $, is also in $D_f$.
\end{remark}

\begin{theorem}\label{thm:null1} Let $\tf \in \tReal[x_1, \dots,
x_n]$, $\tpG_1, \dots ,\tpG_k \in \tTrop[x_1, \dots, x_n]$, and
let $\tIdA $ be the ideal generated by $\tg_1, \dots ,\tg_k$. Then
$\tf \in\sqrt{\tIdA}$ if and only $\cfSet(\tf) \subS
\cfSet(\tIdA)$.
\end{theorem}
Please note that here we work on the reduced tropical semiring
$\tTrop[x_1,\dots,x_n]$, in other word polynomials are identified
with polynomials functions. When the notations are clear from the
context, for short, we write $D \in \cfSet(\tf)$ for a connected
component $D_{\tf} \in \cfSet(\tf)$. (The proof of this theorem
follows after the arguments of \cite[Theorem 3.5 and its
Corollary]{shustinNull}.)

\begin{proof} $(\Rightarrow)$ Assume $\tf \in\sqrt{\tIdA}$, then
$\tf^m = \sum_i \th_i \tg_i$, where $\th_i \in
\tTrop[x_1,\dots,x_n]$, $m \in \Net$. Suppose $D \in
\cfSet(\tf^m)$,  then $ \tf^m |_D$ must coincide with one of the
terms $(\th_i \tg_i)|_D$ in the expression $(\bigoplus_i \th_i
\tg_i) |_D$, since otherwise the the latter function would have a
ghost value inside $D$. Then, by definition, both $\th_i|_D$ and
$\tg_i|_D$ don't have ghost evaluations over $D$, which means $D
\subset D'$ for some $D' \in \cfSet(\tg_i)$. (Recall that
$\cfSet(\tf) = \cfSet(\tf^m)$.)

\medskip
$(\Leftarrow)$ Distribute the connected components of
$\cfSet(\tf)$ into disjoint subsets $\Pi_j$, $j\in J$, where
$J\subset\{1,\dots,k\}$, such that, for any $j\in J$ and
$D\in\Pi_j$, we have $D \subset D'$ for some $D' \in
\cfSet(\tpG_j)$ . Fix some $j\in J$, pick $D \in \Pi_j$, and
assume $\tf|_D = \tf_i |_D$ for some monomial $ \tf_i = \al_i
x^{i_1} \cdots x^{i_n}$. Similarly, we may assume $\tg_j |_{D'} =
\tg_{j,r} |_{D'}$ for some monomial $\tpG_{j,r}=\bt_r x^{r_1}
\cdots x^{r_n}$ of $\tg_j$  (in particular $\tg_j |_{D} =
\tg_{j,r} |_{D}$ where $\tg_{j,r}$ is a tangible monomial).

We claim that for any $t = 1,\dots,n $ \begin{equation}\label{e1}
i_t > 0 \quad \text{ whenever } \quad  r_t > 0;
\end{equation} otherwise, i.e. $i_t =0$ and
$r_t >0$, take a point $\bfa \in D \subset \Trop^{(n)}$  having
only tangible coordinates except the $t$'th coordinate which has a
ghost value (Remark  \ref{rem:con}). Then, $\tpF(\bfa) \in
\etUnit$ on $D$ while $\tg_j(\bfa) \in \Real$ and thus $D
\nsubseteq D'$.

 Condition (\ref{e1}) yields that there is $m_1$
such that, for any $m\ge m_1$ and $D\in\Pi_j$, one has
\begin{equation}\label{eq:n1.5} m \cdot i_t \ge  r_t,  \quad \
t=1, \dots ,n.
\end{equation} Accordingly, we
define the function
\begin{equation}\label{eq:n2} F_{D,m}|_D = \frac{\tf^m | _D}{\tg_j | _D} = \frac{\tf_i^m | _D}{\tg_{j,r} | _D} =
 \frac{\al_i}{\bt_{j,r}}
x_1^{m i_1-r_1} \cdots x_n^{m i_n-r_n} |_D ,\quad D\in\Pi_j,\ m\ge
m_1, \end{equation} for which $m \cdot i_t-r_t$ have always
nonnegative integral values for any $t = 1,\dots,n$. (Note that,
due to \Ref{eq:n1.5}, over $D$ this function is described by a
proper polynomial.)

We claim that there exists $m_2$, such that for any $D \in \Pi_j$,
in the complement of the closure of $D$ (denoted
$\overline{D}^c$), we have
\begin{equation}\label{eq:n3} \tf ^m > F_{D,m} \tg_j\quad\text{whenever}\quad m\ge m_2\ . \end{equation}
Indeed, write $\tf = F  G$, $\tg_j = F' G'$, where $F$, $F'$ are
monomials, and $G = \bigoplus_{\bfk} \gm_{\bfk} \bfx^{\bfk}$, $G'
= \bigoplus_{\bfk'} \gm'_{\bfk'} \bfx^{\bfk'}$ are polynomials
(referred to as functions) equal $0$ along $D$. (In particular, as
functions, $G$ and $G'$ are convex functions.)  Then
$F_{D,m}=\frac{F^m }{F'}$, which  is clearly a monomial on $D$,
and thus $\frac{\tf^m}{F_{D,m} \tg} =\frac{G^m}{G'}$.
By the convexity of $G$, and the fact it equal $0$ on $D$, we have
$G\big|_{\overline {D}^c} > 0$. Since, $G > 0$ and respectively
$k_t \ge 1$, $t=1,...,n$, outside $\overline D$, we obtain
\Ref{eq:n3} when $m_2$ exceeds all the values of the $k'_t$'s of
$\bfk'$ with respect to $G'$.

Define $\th_i=\bigoplus_{D\in\Pi_J}F_{D,m}$. This is a tropical
polynomial as $m\ge m_1$ and, due to Equations \Ref{eq:n2} and
\Ref{eq:n3}, it satisfies
$$(\th_i \tg_j)\big|_D=\tf^m\big|_D,\quad
(\th_i \tg_j)\big|_{\overline{ D}^c}\ <\  \tf^m \big|_{\overline
{D}^c},\quad D\in\Pi_j,\quad m\ge m_2\ ,$$
where $\overline{D}^c$ is the complement of the closure of $D \in
\cfSet(\tf)
 $.
 \end{proof}

\begin{theorem}\label{thm:AlgebraicTropicalNullstellensatz} \textbf{(Algebraic  Nullstellensatz)}
Let $\tIdA \subset \tReal[x_1,\dots,x_n]$ be a finitely generated
tropical ideal, where $\ttUnit[x_1,\dots,x_n] \subseteq \tIdA$,
then
$$\rtIdA = \Idl(\cfSet(\tIdA)).$$
\end{theorem}
\begin{proof} $(\subseteq)$ Assume $\tfF \in
  \rtIdA$, then $\tf^m \in \tIdA$ for some positive $m\in \Net$,
  and hence $\cfSet(\tf^m) \subS
  \cfSet(\tIdA)$. By Lemma \ref{thm:coVarOfPower},
  $\cfSet(\tf)= \cfSet(\tf^m)$ and, since $\cfSet(\tf) \subS
  \cfSet(\tIdA)$, then $\tf \in \Idl(\cfSet(\tIdA))$.

 $(\supseteq)$ When $\tf \in \Idl(\cfSet(\tIdA))$ it
  means that $\cfSet(\tf) \subS
  \cfSet(\tIdA)$, namely,
  each $\conC_{\tf}  \in \cfSet(\tf)$ is contained
  in some component $\conC_{\tIdA}  \in \cfSet(\tIdA)$ and hence in
  some component $\conC_{\tf_i}  \in \cfSet(\tf_i)$ of some $\tf_i \in
  \tIdA$. The proof is then completed by applying Theorem
  \ref{thm:null1}.
\end{proof}


\end{document}